\documentclass[10pt]{article}
\usepackage{amsmath,amsthm,amsfonts,amssymb,enumerate,amscd,graphicx,ytableau,subcaption,mathrsfs,tikz,tikz-cd,float,lmodern}
\usepackage{authblk}
\newcounter{mainthm}

\theoremstyle{plain}
\newtheorem{MTheorem}[mainthm]{Main Result}
\newtheorem*{MTheorem*}{Main Result}
 
\newtheorem{Thm}[subsection]{Theorem}
\newtheorem{Prop}[subsection]{Proposition}

\newtheorem{Cor}[subsection]{Corollary}

\newtheorem{Lem}[subsection]{Lemma}

\theoremstyle{definition}

\newtheorem{defn}[subsection]{Definition}

\newtheorem{Ex}[subsection]{Example}

\newcommand{\La}{L_{\mathsf{A}}(k,N-k)}
\newcommand{\Lpart}{L_{\mathsf{A}}^\text{part}(k,N-k)}
\newcommand{\Ldiag}{L_{\mathsf{A}}^\text{diag}(k,N-k)}
\newcommand{\Lcirc}{L_{\mathsf{A}}^\text{circ}(k,N-k)}
\newcommand{\Ltab}{L_{\mathsf{A}}^\text{tab}(k,N-k)}

\renewcommand{\L}{L_{\mathsf{A}}}
\newcommand{\Lp}{L_{\mathsf{A}}^\text{part}}
\newcommand{\Lt}{L_{\mathsf{A}}^\text{tab}}
\newcommand{\Lc}{L_{\mathsf{A}}^\text{circ}}
\newcommand{\Ld}{L_{\mathsf{A}}^\text{diag}}

\newcommand{\Da}{D_{\mathsf{A}}(k,N-k)}
\newcommand{\Dpart}{D_{\mathsf{A}}^\text{part}(k,N-k)}
\newcommand{\Ddiag}{D_{\mathsf{A}}^\text{diag}(k,N-k)}
\newcommand{\Dcirc}{D_{\mathsf{A}}^\text{circ}(k,N-k)}
\newcommand{\Dtab}{D_{\mathsf{A}}^\text{tab}(k,N-k)}

\newcommand{\D}{D_{\mathsf{A}}}
\newcommand{\Dp}{D_{\mathsf{A}}^\text{part}}
\newcommand{\Dt}{D_{\mathsf{A}}^\text{tab}}
\newcommand{\Dc}{D_{\mathsf{A}}^\text{circ}}
\newcommand{\Dd}{D_{\mathsf{A}}^\text{diag}}

\newcommand{\gamdp}{\gamma^{d\rightarrow p}}
\newcommand{\gampd}{\gamma^{p\rightarrow d}}
\newcommand{\gampt}{\gamma^{p\rightarrow t}}
\newcommand{\gamtp}{\gamma^{t\rightarrow p}}
\newcommand{\gamtc}{\gamma^{t\rightarrow c}}
\newcommand{\gamct}{\gamma^{c\rightarrow t}}
\newcommand{\m}{\mathsf{m}}

\makeatletter
\newcommand\setItemnumber[1]{\setcounter{enum\romannumeral\@enumdepth}{\numexpr#1-1\relax}}
\makeatother

\newcommand{\ps}{\mathbf{s}}
\newcommand{\pt}{\mathbf{t}}
\newcommand{\ds}{\mathsf{s}}
\newcommand{\dt}{\mathsf{t}}

\newcommand{\gs}{\sigma}
\newcommand{\gt}{\tau}

\newcommand{\BZ}{\mathbb{Z}}

\newcommand{\myA}{\mathsf{A}}
\newcommand{\myB}{\mathsf{B}}
\newcommand{\myC}{\mathsf{C}}
\newcommand{\myD}{\mathsf{D}}
\newcommand{\relt}{\mathbf{r}}
\newcommand{\selt}{\mathbf{s}} 
\newcommand{\telt}{\mathbf{t}}
\newcommand{\uelt}{\mathbf{u}} 
\newcommand{\velt}{\mathbf{v}}
\newcommand{\welt}{\mathbf{w}} 
\newcommand{\xelt}{\mathbf{x}}
\newcommand{\yelt}{\mathbf{y}} 
\newcommand{\zelt}{\mathbf{z}}

\newcommand{\subsc}[1]{\textnormal{\tiny \textsc{#1}}}
\definecolor{blue}{RGB}{66,134,244}
\definecolor{green}{RGB}{37,150,41}
\definecolor{orange}{RGB}{247,169,0}
\definecolor{purple}{RGB}{152,0,247}
\definecolor{gold}{RGB}{244,225,12}
\definecolor{magenta}{RGB}{244,12,206}
\definecolor{cyan}{RGB}{12, 233, 244}
\definecolor{brick}{RGB}{99,27,8}
\definecolor{seafoam}{RGB}{152,234,152}
\definecolor{lavender}{RGB}{255,175,255}

\newcommand{\tab}[2]{
\begin{ytableau}
#1\\
#2
\end{ytableau}
}

\newcommand{\partpic}[1]{
\ydiagram[*(gray)]
{#1}
*[*(white)]{4,4}}

\newcommand{\circdiag}[1]{
\ydiagram[*(white)\bullet]
{#1}
*[*(white)]{3,3}}

\newcommand{\TikzPartTwoThree}[1]{
\ydiagram[*(white) \bullet]
{#1}
*[*(white)]{3,3}}

\usetikzlibrary{matrix}

\newenvironment{enumeratei}{\begin{enumerate}[\upshape (i)]}
        {\end{enumerate}}

\numberwithin{equation}{section}
\title{Diamond-colored distributive lattices, move-minimizing games, and fundamental Weyl symmetric functions: The type $\myA$ case} 
\author{Robert G. Donnelly}
\author{Elizabeth A. Donovan}
\author{Timothy A. Schroeder}
\affil{Department of Mathematics \& Statistics\\
Murray State University\\
Murray, KY  42071\quad USA
\vskip.5cm
\small \tt rdonnelly@murraystate.edu\quad edonovan@murraystate.edu\quad tschroeder@murraystate.edu}
\date{\today}

\begin{document}

\maketitle

\begin{abstract}    
\sloppy We present some elementary but foundational results concerning diamond-colored modular and distributive lattices and connect these structures to certain one-player combinatorial ``move-minimizing games,'' in particular, a so-called ``domino game.''  The objective of this game is to find, if possible, the least number of ``domino moves'' to get from one partition to another, where a domino move is, with one exception, the addition or removal of a domino-shaped pair of tiles.  We solve this domino game by demonstrating the somewhat surprising fact that the associated ``game graphs'' coincide with a well-known family of diamond-colored distributive lattices which shall be referred to as the ``type $\myA$ fundamental lattices.'' 
These lattices arise as supporting graphs for the fundamental representations of the special linear Lie algebras and as splitting posets for type $\myA$ fundamental symmetric functions, connections which are further explored in sequel papers for types $\myA$, $\myC$, and $\myB$. 
In this paper, this connection affords a solution to the proposed domino game as well as new descriptions of the type $\myA$ fundamental lattices.
\end{abstract}
%


\ytableausetup{smalltableaux}
\section{Introduction}\label{s:intro}

Consider this well-known logic puzzle:  Given an $8\times 8$ black and white checkerboard with 64 squares, let a domino tile be a rectangle formed by taking two checkerboard-size squares and joining them along an edge.  Is it possible to perfectly pave the checkerboard with Domino tiles?  The answer to this question is, of course, ``Yes.'' 

But if two opposite corners from the checkerboard were removed to yield a new board with 62 squares, would it be possible to perfectly tile this modified board?  The answer to this is ``No.''  Since each domino tile placed on the board covers a red square and a white square on the checkerboard, any tiling must cover the same number of red and white squares.  But opposite corners of the board are the same color, so the modified board will not have the same number of red and white squares.  

This simple logic puzzle illustrates a couple of ideas that will be useful in what follows:  (1) We will propose a game that uses domino tiles and pavings of certain shapes within a certain rectangular grid of squares, and (2) Although checkerboard-style coloring is not needed to state the possiblility/impossibility problems we address, it will be useful in formulating our solution.

To arrive at our game, we fix positive integers $k$ and $N$ with $1\leq k\leq N-1$.  We will define a domino-type game using tiles like those described above, however, our game board will consist of a $k\times (N-k)$ grid of squares.  The Domino Game is played with \emph{singletons} (marked tiles each the size of a game board square) that are to be added/removed in \emph{domino pairs} (two square tiles meeting along an edge) or, in certain limited situations, one at a time.  

A \emph{legal shape} is a paving of all or part of our game board with singletons such that 
\begin{enumerate}
	\item The paved squares on a given row of the game board are \emph{left-justified}, and
	\item The number of paved squares on a given row is \emph{at least as large as} the number of paved squares on the next row.
\end{enumerate}
All of this is to say that each each legal shape is the ``partition diagram'' (or ``Ferrer's diagram'') for a partition with no more than $k$ parts, and largest part of size at most $N-k$.  We call this partition, or its corresponding shape, a \emph{$k\times(N-k)$ partition}.   

A \emph{domino move} is the addition or removal of 
\begin{itemize}
	\item A domino pair, or
	\item A singleton in the upper right corner only.
\end{itemize}
When applied to a legal shape on a given game board, a domino move is \emph{legal} if it results in another legal shape for the same board.

\begin{figure}[h]
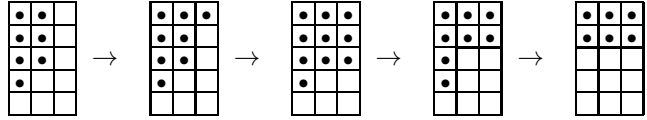

\centering
\ytableausetup{centertableaux}

\ydiagram[*(white) \bullet]
{2,2,2,1}
*[*(white)]{3,3,3,3,3}
\hspace{.1cm}
$\rightarrow$
\hspace{.2cm}
\ydiagram[*(white) \bullet]
{3,2,2,1}
*[*(white)]{3,3,3,3,3}
\hspace{.1cm}
$\rightarrow$
\hspace{.2cm}
\ydiagram[*(white) \bullet]
{3,3,3,1}
*[*(white)]{3,3,3,3,3}
\hspace{.1cm}
$\rightarrow$
\hspace{.2cm}
\ydiagram[*(white) \bullet]
{3,3,1,1}
*[*(white)]{3,3,3,3,3}
\hspace{.1cm}
$\rightarrow$
\hspace{.2cm}
\ydiagram[*(white) \bullet]
{3,3}
*[*(white)]{3,3,3,3,3}
\caption{A game consisting of four legal moves to move from the shape on the left to the shape on the right.}
\label{fig:examplegame}
\end{figure}

The \emph{Domino Game} begins with a game board of fixed size and two legal shapes.  The goal of the game is to get from one shape to the other using only \emph{legal domino moves}.  See Figure \ref{fig:examplegame} for an example of an instance of a successful domino game played on a $5\times (8-5)$ game board.

The Domino Game gives rise to three natural questions:
\begin{itemize}
	\item Is it always possible to get from one legal shape to another using legal domino moves?
	\item If so, what is the smallest number of moves required?
	\item And if so, can a smallest sequence of moves be explicitly prescribed?
\end{itemize}
To answer these questions is to solve the Domino Game. 

This Domino Game will serve as our primary example of a ``move-minimizing game''. However, at the outset of the paper in Section \ref{s:lattices}, we will concern ourselves with the technical and seemingly ancillary topic of diamond-colored modular and distributive lattices. These arise naturally in certain areas of algebraic combinatorics: as supporting graphs, which are combinatorial models for complex semi-simple Lie algebra representations (see for example \cite{Don2}), and as splitting posets, which are poset models for so-called Weyl symmetric functions (see for example \cite{ADLMPPW}). 

In Section \ref{s:classicalL}, we present a classical family of diamond-colored distributive lattices, which we call here the \emph{type $\myA$ fundamental lattices.}  We also provide a new description of these lattices using what we call \emph{diagonal coordinates.} 

In Section \ref{s:games}, we describe the general notion of a move-minimizing game, and then explore a solution to any such game whose associated graph is the order diagram for a diamond-colored modular or distributive lattice. We apply this methodology in Section \ref{s:dominos} to solve our proposed Domino Game by realizing the underlying game graph as a type $\myA$ fundamental lattice. 

Although this paper is almost exclusively combinatorial, and the first part highly technical, it is part of a series of papers whose aim is also algebraic; namely to present new combinatorial models for the fundamental representations of the classical Lie algebras: the type $\myA$ or special linear Lie algebras; the type $\myB$ or odd orthogonal Lie algebras; the type $\myC$ or symplectic Lie algebras; and the type $\myD$ or even orthogonal Lie algebras, and their associated Weyl symmetric functions. 

The purpose of this paper is to provide the necessary foundational results concerning diamond-colored modular and distributive lattices and showcase our perspective with the presentation of the type $\myA$ fundamental lattices via various natural coordinatizations as well as via the Domino Game. The paper will have two main results.  
\begin{MTheorem}\label{t:intromain1} Let $\Gamma$ be a diamond-colored distributive lattice.  Let $\ps$ and $\pt$ be elements of $\Gamma$.  Then 
\begin{enumeratei}
	\item A shortest path from $\ps$ to $\pt$ can be chosen to pass through $\ps\wedge\pt$ (or $\ps\vee\pt$).
	\item The length of such a path can be calculated in terms of $\ps$ and $\pt$.
	\item The path can be explicitly prescribed in terms of $\ps$ and $\pt$.  
\end{enumeratei}
\end{MTheorem}
\begin{MTheorem}\label{t:intromain2} The game graph associated to the Domino Game is a diamond-colored distributive lattice.    
\end{MTheorem}
Result \ref{t:intromain1} is actually a culmination of Theorems \ref{t:FundamentalTheorem}, \ref{t:FullLengthTheorem} and \ref{t:MoveMinGameTheorem} from Sections \ref{s:lattices} and \ref{s:games}.  Result \ref{t:intromain2} is shown in Section \ref{s:dominos} as Theorem \ref{t:main}, and will be the basis to solving the Domino game via the methodology put forward in the first portion of the paper.

%
%

The combinatorial structures considered in this paper and in the sequels \cite{dds2} and \cite{bdds} are finite posets which we view as graphs, when identified with their order diagrams.  Most often, such a graph will either have its edges or its vertices colored by elements from some index set $I$ (usually a set of positive integers).  The conventions and notation we use concerning such structures largely borrow from \cite{Don2}, \cite{ADLP1}, \cite{ADLMPPW}, and \cite{Stan}.  The reader should note that in this paper, many proofs are omitted.  The details of such proofs are technical, but overall they are in fact quite intuitive and require only elementary reasoning.

\section{Diamond-colored modular \& distributive lattices.}\label{s:lattices} 

In this section we establish some foundational results for our study of ``diamond-colored'' modular and distributive lattices.  The purpose of this section is two-fold:  (1) We will review some classic terminology related to edge- and vertex-colored graphs and connect these ideas to new concepts related to distributive lattices; and (2) The results in this section will provide the theoretical foundation for our later investigation of ``move-minimizing'' games in general and the domino game in particular.   This section is by nature quite technical, and the reader interested solely in the game questions put forward in Section \ref{s:intro} is encouraged to move ahead to Sections \ref{s:games} and \ref{s:dominos}.

\paragraph{Graphs and posets.}  Suppose $\Gamma$ is a simple directed graph with vertex set $\mathcal{V}(\Gamma)$ and directed edge set $\mathcal{E}(\Gamma)$.  If $\Gamma$ is accompanied by a function $EC_\Gamma:\mathcal{E}(\Gamma) \longrightarrow I$ (respectively $VC_{\Gamma}: \mathcal{V}(\Gamma) \longrightarrow I$), then we say $\Gamma$ is {\em edge-colored}  (resp.\ {\em vertex-colored}) {\em by the set} $I$. 
An isomorphism of two edge-colored (respectively, vertex-colored) directed graphs is a bijection $\phi$ of their vertex sets such that $\phi$ and $\phi^{-1}$ both preserve directed edges and edge colors (resp., vertex colors). 

If $\Gamma$ arises as the order diagram for a poset $R$ with partial ordering ``$\leq$'', then we identify $\Gamma$ with $R$ and refer to $R$ as an edge-colored (resp.\ vertex-colored) poset. In such an edge-colored poset $R$, by definition we have $\xelt \xrightarrow{i} \yelt$ for some vertices $\xelt$ and $\yelt$ and some color $i \in I$ if and only if (a) $\xelt < \yelt$ with respect to the partial ordering on $\mathcal{V}(R)$ and (b) for any other vertex $\zelt$ we have $\xelt = \zelt$ or $\zelt = \yelt$ whenever $\xelt \leq \zelt \leq \yelt$.  For this $i$-colored edge $\xelt \xrightarrow{i} \yelt$, we say that $\yelt$ is {\em above} $\xelt$ and that $\xelt$ is \emph{below} $\yelt$. In figures, when such an edge is depicted without an arrowhead, the implied direction is ``up.'' 

\sloppy\paragraph{Mountain and valley paths.}  Regard $R$ to be a poset that is edge-colored by a set $I$.  A {\em path} from $\selt$ to $\telt$ in $R$ is a sequence $\mathcal{P} = (\selt = \xelt_{0}, \xelt_{1}, \ldots, \xelt_{k} = \telt)$ such that for $1 \leq j \leq k$ we have either $\xelt_{j-1} \xrightarrow{i_j} \xelt_{j}$ or $\xelt_{j} \xrightarrow{i_j} \xelt_{j-1}$, where $(i_{j})_{j=1}^{k}$ is a sequence of colors from $I$.  
This path has {\em length} $k$, written $\mathbf{\ell}(\mathcal{P})$, and we allow paths to have length $0$. For any $i \in I$, we let $a_{i}(\mathcal{P}) := \left|\{j \in [k]\, \mid \, \xelt_{j-1} \xrightarrow{i} \xelt_{j} \mbox{ in $\mathcal{P}$} \}\right|$, a count of ``ascending'' edges of color $i$ in the path, and $d_{i}(\mathcal{P}) := \left|\{j \in [k], \mid \, \xelt_{j} \xrightarrow{i} \xelt_{j-1} \mbox{ in $\mathcal{P}$}\}\right|$, a count of ``descending'' edges of color $i$, where $[k] :=\{1, 2, \ldots, k\}$. 
Of course, $\mathbf{\ell}(\mathcal{P}) = \sum_{i \in I}(a_{i}(\mathcal{P}) + d_{i}(\mathcal{P}))$. Say $\mathcal{P}$ is {\em simple} if each vertex appearing in the path appears exactly once. 

If $\selt$ and $\telt$ are within the same connected component of $R$, then the distance $d(\selt,\telt)$ between $\selt$ and $\telt$ is the minimum length achieved when all paths from $\selt$ to $\telt$ in $R$ are considered; any minimum-length-achieving path is a {\em shortest path}.  Our poset $R$ has the {\em diamond coloring property} if whenever  
\parbox{1.1cm}{\begin{center}
\setlength{\unitlength}{0.2cm}
\begin{picture}(4,3.5)
\put(2,0){\circle*{0.5}} \put(0,2){\circle*{0.5}}
\put(2,4){\circle*{0.5}} \put(4,2){\circle*{0.5}}
\put(0,2){\line(1,1){2}} \put(2,0){\line(-1,1){2}}
\put(4,2){\line(-1,1){2}} \put(2,0){\line(1,1){2}}
\put(0.15,0){\em \small k} \put(3.1,0){\em \small l}
\put(.2,3){\em \small i} \put(3.1,3){\em \small j}
\end{picture} \end{center}} is an edge-colored subgraph of the order  
diagram for $R$, then $i = l$ and $j = k$.  

A simple path $\mathcal{P}$ in $R$ is a {\em mountain path}  from $\selt$ to $\telt$ if for some $j \in \{0,1,\dots,k\}$ and $\uelt := \xelt_{j}$ we have $\selt=\xelt_{0} \rightarrow \xelt_{1} \rightarrow \cdots \rightarrow \uelt \leftarrow \cdots \leftarrow \xelt_{k} = \telt$, in which case we call $\uelt$ the {\em apex} of the mountain path. 

Similarly, the simple path $\mathcal{P}$ is a {\em valley path} from $\selt$ to $\telt$ if for some $j \in \{0,1,\dots,k\}$ and $\velt := \xelt_{j}$ we have $\selt=\xelt_{0} \leftarrow \xelt_{1} \leftarrow \cdots \leftarrow \velt \rightarrow \cdots \rightarrow \xelt_{k} = \telt$, in which case we call $\velt$ the {\em nadir} of the mountain path. 

Our poset $R$ is {\em topographically balanced} if (1) whenever $\velt \rightarrow \selt$ and $\velt \rightarrow \telt$ for distinct $\selt$ and $\telt$ in $R$, 
then there exists a unique $\uelt$ in $R$ such that $\selt \rightarrow \uelt$ and $\telt \rightarrow \uelt$, and (2) whenever  $\selt \rightarrow \uelt$ and $\telt \rightarrow \uelt$ for distinct $\selt$ and $\telt$ in $R$, then there exists a unique $\velt$ in $R$ such that $\velt \rightarrow \selt$ and $\velt \rightarrow \telt$. Informally, this just says that any length two mountain path that is not a chain is uniquely balanced by a length two valley path that is not a chain, and vice-versa. If $R$ is topographically balanced, any simple path $\mathcal{P} = (\selt = \xelt_{0}, \xelt_{1}, \ldots, \xelt_{k} = \telt)$ from $\selt$ to $\telt$ can be {\em mountain-ized} to form a new path $\mathcal{P}_{\mbox{\em \scriptsize mountain}}$ as follows: 
\begin{enumerate}
\item If there is no $j \in [k-1]$ such that $\xelt_{j-1} \leftarrow \xelt_{j} \rightarrow \xelt_{j+1}$, then return $\mathcal{P}$ as $\mathcal{P}_{\mbox{\em \scriptsize mountain}}$. 
\item Otherwise, let $j \in [k-1]$ be least such that $\xelt_{j-1} \leftarrow \xelt_{j} \rightarrow \xelt_{j+1}$ and form a new path $\mathcal{P}'$ be replacing $\xelt_{j}$ with the unique $\xelt'_{j}$ for which $\xelt_{j-1} \rightarrow \xelt'_{j} \leftarrow \xelt_{j+1}$. 
\item Return to the first step of the process, using $\mathcal{P}'$ as $\mathcal{P}$. 
\end{enumerate}
The mountain-ization $\mathcal{P}_{\mbox{\em \scriptsize mountain}}$ has the salient properties that it is a mountain path and that $\mathbf{\ell}(\mathcal{P}) = \mathbf{\ell}(\mathcal{P}_{\mbox{\em \scriptsize mountain}})$. 
If  $R$ is also diamond-colored by a set $I$, then $a_{i}(\mathcal{P}) = a_{i}(\mathcal{P}_{\mbox{\em \scriptsize mountain}})$ and $d_{i}(\mathcal{P}) = d_{i}(\mathcal{P}_{\mbox{\em \scriptsize mountain}})$ for each $i \in I$. 
Similarly define the {\em valley-ization} $\mathcal{P}_{\mbox{\em \scriptsize valley}}$ of $\mathcal{P}$. 

\paragraph{Results.}  We now connect the ideas of moutain- and valley-ization of paths with the classical concepts of \emph{ranked lattices}, in particular \emph{distributive} and \emph{modular lattices}.  (For reference, see \cite{davey}.)  In our setting, if $\rho$ denotes a rank function for some poset $R$, then for any path $\mathcal{P}$ from $\selt$ to $\telt$, we have \[\rho(\telt) - \rho(\selt) =  \sum_{i \in I}\big(a_{i}(\mathcal{P}) - d_{i}(\mathcal{P})\big),\]  
which serves as an expression for the rank of $\telt$ whenever $\rho(\selt) = 0$.

\begin{Thm}\label{t:topbalancerank} A lattice $L$ is topographically balanced if and only if $L$ is ranked with unique rank function $\rho$ satisfying 
\[2\rho(\selt \vee \telt) - \rho(\selt) - \rho(\telt) = \rho(\selt) + \rho(\telt) - 2\rho(\selt \wedge \telt)\] 
for all $\selt, \telt \in L$.
\end{Thm}

\begin{Thm}\label{t:ModularLatticeTheorem} 
Let $\selt$ and $\telt$ be elements of a topographically balanced lattice $L$ and let $\mathcal{P}$ be any simple path in $L$ from $\selt$ to $\telt$.  Then the following are equivalent:
\begin{enumerate}
	\item $\mathcal{P}$ is a path of shortest length from $\selt$ to $\telt$.
	\item $\displaystyle \mathbf{\ell}(\mathcal{P}) = 2\rho(\selt \vee \telt) - \rho(\selt) - \rho(\telt) = \rho(\selt) + \rho(\telt) - 2\rho(\selt \wedge \telt)$.
	\item $\mathcal{P}_{\mbox{\em \scriptsize mountain}}$ has apex $\uelt = \selt \vee \telt$ and $\mathcal{P}_{\mbox{\em \scriptsize valley}}$ has nadir $\velt = \selt \wedge \telt$. 
\end{enumerate}
\end{Thm}


The fact that a distributive lattice is ranked is a consequence of the following standard result, whose proof amounts to the observation that a distributive lattice is topographically balanced. 

\begin{Lem}\label{l:DistributiveIsModularLemma} Any distributive lattice is modular.
\end{Lem}

The following Lemma connects modular and diamond-colored lattices.  It is a utility that aids in proofs of other results (e.g. Theorem \ref{t:FullLengthTheorem}) and helps in analyzing movement between vertices in diamond-colored modular lattices (e.g.\ computing the rank of an element). 

\begin{Lem}\label{l:ColorsLemma} Let $L$ be a modular lattice diamond-colored by a set $I$.  Suppose $\selt \leq \telt$.  Suppose $\mathcal{P} = (\selt = \relt_{0} \xrightarrow{i_{1}} \relt_{1} \xrightarrow{i_{2}} \relt_{2} \xrightarrow{i_{3}} \cdots \xrightarrow{i_{p-1}} \relt_{p-1} \xrightarrow{i_{p}} \relt_{p} = \telt)$ and  $\mathcal{Q} = (\selt = \relt'_{0} \xrightarrow{j_{1}} \relt'_{1} \xrightarrow{j_{2}} \relt'_{2} \xrightarrow{j_{3}} \cdots \xrightarrow{j_{q-1}} \relt'_{q-1} \xrightarrow{j_{q}} \relt'_{q} = \telt)$ are two paths from $\selt$ up to $\telt$. Then, $p = q$, and $a_{i}(\mathcal{P}) = a_{i}(\mathcal{Q})$ for all $i \in I$. Moreover, if $\relt_{1}$ and $\relt'_{p-1}$ are incomparable, then $i_{1} = j_{p}$.
\end{Lem} 

\paragraph{Diamond-colored distributive lattices and order ideals.} The following discussion of diamond-colored distributive lattices and certain related vertex-colored posets is quite technical, but it does encompass the classical uncolored situation (see for example Section 3.4 of \cite{Stan}). These concepts have antecedents in work of Proctor and Stembridge (see e.g. \cite{Proc2}, \cite{proc3}, \cite{StemQ}, \cite{Stem}), but there seems to be no standard treatment of these ideas.  

A diamond-colored distributive lattice can be constructed as follows. Let $P$ be a poset with vertices colored by a set $I$.  An {\em order ideal} $\xelt$ from $P$ is a vertex subset of $P$ with the property that $u \in \xelt$ whenever $v \in \xelt$ and $u \leq v$ in $P$.  Let $L$ be the set of order ideals from $P$.  For $\xelt, \yelt \in L$, write $\xelt \leq \yelt$ if and only if $\xelt \subseteq \yelt$.  
  
With respect to this partial ordering, $L$ is a distributive lattice: $\xelt \vee \yelt = \xelt \cup \yelt$ and $\xelt \wedge \yelt = \xelt \cap \yelt$ for all $\xelt, \yelt \in L$.  One can easily see that $\xelt \rightarrow \yelt$ in $L$ if and only if $\xelt \subset \yelt$ and $\yelt - \xelt = \{v\}$ for some maximal element $v$ of $\yelt$, with $\yelt$ thought of as a subposet of $P$ in the induced order.  

In this case, we declare that $EC_{L}(\xelt \rightarrow \yelt) := VC_{P}(v)$, thus making $L$ an edge-colored distributive lattice: viewed as an element of $\mathbf{J}_{color}(P)$, an order ideal $\xelt$ from $P$ is above an edge of color $i$ if and only if $\xelt$ has a maximal (resp.\ minimal) element of color $i$. One can easily check that $L$ has the diamond-coloring property.  The diamond-colored distributive lattice just constructed is given special notation: we write $L := \mathbf{J}_{color}(P)$.  
Note that if $P \cong Q$ as vertex-colored posets, then $\mathbf{J}_{color}(P) \cong \mathbf{J}_{color}(Q)$ as edge-colored posets.  Moreover, $L$ is ranked with rank function given by $\rho(\telt) = |\telt|$.  In particular, the length of $L$ is $|P|$. 


The process described in the previous paragraph can be reversed. Given a diamond-colored distributive lattice $L$, an element $\xelt$ is {\em join irreducible} if $\xelt \not= \min(L)$ and whenever $\xelt = \yelt \vee \zelt$ then $\xelt = \yelt$ or $\xelt = \zelt$.  One can see that $\xelt$ is join irreducible if and only if $\xelt$ 
covers precisely one other vertex in $L$, i.e.\ $\left|\{\xelt'\, \mid\, \xelt' \rightarrow \xelt\}\right| = 1$.  Let $P$ be the set of all join irreducible elements of $L$ with the induced partial ordering.  Color the vertices of $P$ by the rule: $VC_{P}(\xelt) := EC_{L}(\xelt' \rightarrow \xelt)$. We call $P$ the {\em vertex-colored poset of join irreducibles} and denote it by $P := \mathbf{j}_{color}(L)$.  If $K \cong L$ is an isomorphism of diamond-colored lattices, then $\mathbf{j}_{color}(K) \cong \mathbf{j}_{color}(L)$ is an isomorphism of vertex-colored posets.  


What follows is a dual to the above constructions of diamond-colored distributive lattices. A {\em filter} from a vertex-colored poset $P$ is a subset $\xelt$ with the property that if $u \in \xelt$ and $u \leq v$ in $P$ then $v \in \xelt$.  Note that for $\xelt \subseteq P$, $\xelt$ is a filter if and only if the set complement $P - \xelt$ is an order ideal.  Now partially order all filters from $P$ by reverse containment: $\xelt \leq \yelt$ if and only if $\xelt \supseteq \yelt$ for filters $\xelt, \yelt$ 
from $P$.  The resulting partially ordered set $L$ is a distributive lattice.  Color the edges of $L$ as in the case of order ideals.  That is, $\xelt$ is below an edge of color $i$ if and only if $\xelt$ has a minimal element of color $i$. The result is a diamond-colored distributive lattice which we denote by $L  = \mathbf{M}_{color}(P)$.  Further, for any filter $\xelt$ from $P$, the rank of $\xelt$ in $\mathbf{M}_{color}(P)$ is $l-|\xelt|$. 

In the other direction, given a diamond-colored distributive lattice $L$, we say $\xelt \in L$ is {\em meet irreducible} if and only if $\xelt \not= \max(L)$ and 
whenever $\xelt = \yelt \wedge \zelt$ then $\xelt = \yelt$ or $\xelt = \zelt$.  One can see that $\xelt$ is meet irreducible if and only if $\xelt$ is covered by exactly one other vertex in $L$.  Now consider the set $P$ of meet irreducible elements in $L$ with the order induced from $L$.  Color the vertices of $P$ in the same way we colored the vertices of the poset of join irreducibles.  The vertex-colored poset $P$ is the {\em poset of meet irreducibles} for $L$.  In this case, we write $P = \mathbf{m}_{color}(L)$.  We have $\mathbf{M}_{color}(P) \cong \mathbf{M}_{color}(Q)$ if $P$ and $Q$ are isomorphic  vertex-colored posets.  We also have $\mathbf{m}_{color}(L) \cong \mathbf{m}_{color}(K)$ if $L$ and $K$ are isomorphic diamond-colored distributive lattices.  

We encapsulate and modestly extend the discussion of the preceding paragraphs in the following  two-part proposition. 

\begin{Prop}\label{p:EncapsulateProposition}  \textbf{(1)} Let $L$ be a distributive lattice diamond-colored by a set $I$.  Then each of $\mathbf{j}_{color}(L)$ and $\mathbf{m}_{color}(L)$ are vertex-colored by the set $I$, and if $L$ has length $l$ with respect to its rank function, then each of these posets has $l$ elements.  \\
\textbf{(2)} Let $P$ be poset with vertices colored by the set $I$.  Then $\mathbf{J}_{color}(P)$ and $\mathbf{M}_{color}(P)$ are distributive lattices that are diamond-colored by the set $I$, and if $P$ has $l$ elements, then each of $\mathbf{J}_{color}(P)$ and $\mathbf{M}_{color}(P)$ has length $l$ with respect to its rank function.  Finally, the unique maximal element of $\mathbf{J}_{color}(P)$ (respectively $\mathbf{M}_{color}(P)$) is $P$ (resp.\ $\emptyset$) and the unique minimal element is $\emptyset$ (resp.\ $P$).
\end{Prop}


The following Theorem shows that the operations $\mathbf{J}_{color}$ (respectively, $\mathbf{M}_{color}$) and $\mathbf{j}_{color}$ (respectively, $\mathbf{m}_{color}$)  are inverses in a certain sense. This is a straightforward generalization of the fundamental theorem of finite distributive lattices. 

\begin{Thm}\label{t:FundamentalTheorem} \textup{(Fundamental Theorem of Finite Diamond-colored Distributive 
Lattices)}
\begin{enumerate}
	\item Let  $L$ be a diamond-colored distributive lattice. Then   
\[L \cong \mathbf{J}_{color}(\mathbf{j}_{color}(L)) \cong 
\mathbf{M}_{color}(\mathbf{m}_{color}(L)).\]
\item Let $P$ be a  
vertex-colored poset.  Then
\[P \cong 
\mathbf{j}_{color}(\mathbf{J}_{color}(P)) \cong 
\mathbf{m}_{color}(\mathbf{M}_{color}(P)).\] 
\end{enumerate}

\end{Thm}

\begin{Cor}\label{c:FirstCorollary} An edge-colored distributive lattice $L$ is isomorphic to $\mathbf{J}_{color}(P)$ or $\mathbf{M}_{color}(P)$ for some vertex-colored poset $P$ if and only if $L$ is diamond-colored.
\end{Cor}


The next corollary states for the record how $\mathbf{J}_{color}$, $\mathbf{j}_{color}$, $\mathbf{M}_{color}$, and $\mathbf{m}_{color}$ interact with the standard vertex- and edge-colored poset operations $*$ (dual), $\sigma$ (recoloring of vertices or edges), $\oplus$ (disjoint union), and $\times$ (Cartesian product). 

\begin{Cor}\label{c:JMCorollary} Let $P$ and $Q$ be posets with vertices colored by a set $I$, and let $K$ and $L$ be diamond-colored distributive lattices with edges colored by $I$. In what follows, $*$, $\sigma$, $\oplus$, $\times$, and $\cong$ account for colors on vertices/edges as appropriate.
\begin{enumerate}
	\item If $K \cong L$, then $\mathbf{j}_{color}(K) \cong \mathbf{m}_{color}(K) \cong \mathbf{m}_{color}(L) \cong \mathbf{j}_{color}(L)$.
	\item If $P \cong Q$, then $\mathbf{J}_{color}(P) \cong \mathbf{M}_{color}(P) \cong \mathbf{M}_{color}(Q) \cong \mathbf{J}_{color}(Q)$.
	\item $\mathbf{J}_{color}(P^{*}) \cong (\mathbf{J}_{color}(P))^{*}$, $\mathbf{J}_{color}(P^{\sigma}) \cong (\mathbf{J}_{color}(P))^{\sigma}$, and $\mathbf{J}_{color}(P \oplus Q) \cong \mathbf{J}_{color}(P) \times \mathbf{J}_{color}(Q)$.
\item $\mathbf{M}_{color}(P^{*}) \cong (\mathbf{M}_{color}(P))^{*}$, $\mathbf{M}_{color}(P^{\sigma}) \cong  (\mathbf{M}_{color}(P))^{\sigma}$, and $\mathbf{M}_{color}(P \oplus Q) \cong \mathbf{M}_{color}(P) \times \mathbf{M}_{color}(Q)$.
\item $\mathbf{j}_{color}(L^{*}) \cong (\mathbf{j}_{color}(L))^{*}$, $\mathbf{j}_{color}(L^{\sigma}) \cong (\mathbf{j}_{color}(L))^{\sigma}$, and $\mathbf{j}_{color}(L \times K) \cong \mathbf{j}_{color}(L) \oplus \mathbf{j}_{color}(K)$. 
\item $\mathbf{m}_{color}(L^{*}) \cong (\mathbf{m}_{color}(L))^{*}$, $\mathbf{m}_{color}(L^{\sigma}) \cong (\mathbf{m}_{color}(L))^{\sigma}$, and $\mathbf{m}_{color}(L \times K) \cong \mathbf{m}_{color}(L) \oplus \mathbf{m}_{color}(K)$.
\end{enumerate}
\end{Cor} 

\paragraph{Subposets and sublattices.} The results that close this section require certain notions of substructures. Assume $Q$ is a weak subposet of $R$.  If, in addition, $Q$ and $R$ are vertex-colored (respectively, edge-colored) by a set $I$ and $VC_{Q}^{-1}(i) \subseteq VC_{R}^{-1}(i)$ (resp.\ $EC_{Q}^{-1}(i) \subseteq EC_{R}^{-1}(i)$) for all $i \in I$, then $Q$ is a {\em vertex-colored} (resp.\ {\em edge-colored}) {\em weak subposet}. 

Let $L$ be a lattice with sublattice $K$. If $K$ is an edge-colored weak subposet of $L$, then we call $K$ an {\em edge-colored sublattice} of $L$.  If both are ranked and have the same length, then say $K$ is a {\em full-length sublattice} of $L$. The next Lemma gives us one way to know whether of a sublattice are also edges of the `parent' lattice.

\begin{Lem}\label{l:FullLengthLemma}  Let $K$ be a full-length sublattice of $L$.  Let $\rho^{(K)}$ and $\rho^{(L)}$ denote the rank functions of $K$ and $L$ respectively.  Then $\rho^{(K)}(\xelt) = \rho^{(L)}(\xelt)$ for all $\xelt$ in $K$, and moreover for all $\xelt$ and $\yelt$ in $K$ we have $\xelt \rightarrow \yelt$ in $K$ if and only if $\xelt \rightarrow \yelt$ in $L$.
\end{Lem}
 
Here is a situation in which a full-length sublattice can be discerned.

\begin{Prop}\label{p:FullLengthWithinProduct} Suppose $L_{1}, L_{2}, \ldots , L_{p}$ are all modular (respectively, distributive) lattices that are diamond-colored by a set $I$, with respective rank functions $\rho^{(1)}, \rho^{(2)}, \ldots , \rho^{(p)}$ and lengths $l^{(1)}, l^{(2)}, \ldots , l^{(p)}$.  Let $L := L_{1} \times L_{2} \times \cdots \times L_{p}$ with inherited rank function $\rho$.\\
\begin{enumerate}
	\item\label{p:FullLengthWithinProduct1} Then $L$ is a modular (resp.\ distributive) lattice, is diamond-colored by $I$, and has length given by $\sum_{q=1}^{p}l^{(q)}$.  Moreover, $\max(L)=(\max(L_{1}),\ldots,\max(L_{p}))$ while $\min(L)=(\min(L_{1}),\ldots,\min(L_{p}))$.  For any $\selt = (s_1, s_2, \ldots , s_p) \in L$ we have $\rho(\selt) = \sum_{q=1}^{p}\rho^{(q)}(s_{q})$.  For any other $\telt = (t_1, t_2, \ldots , t_p)$, we have $\selt \vee_{L} \telt = (s_{1} \vee t_{1}, \ldots , s_{p} \vee t_{p})$ and $\selt \wedge_{L} \telt = (s_{1} \wedge t_{1}, \ldots , s_{p} \wedge t_{p})$.
	\item\label{p:FullLengthWithinProduct2} Suppose $K$ is some vertex subset of $L$ which is closed under component-wise joins and meets.  Further, suppose that $\min(L)$ and $\max(L)$ are in $K$ and there is a path in $L$ from $\min(L)$ to $\max(L)$ that uses only vertices from $K$.  Then $K$ is a full-length sublattice of $L$, is modular (resp.\ distributive), and is diamond-colored by the set $I$.
\end{enumerate}
\end{Prop}

In the special case of diamond-colored distributive lattices, the following result, which is a diamond-colored version of Remark 2.1 of \cite{Don2}, can be applied to help find nice presentations of posets of join irreducibles for distributive lattices which arise as full-length sublattices of larger and more easily described distributive lattices (see e.g.\ \cite{Don2}, \cite{gil}).

\begin{Thm}\label{t:FullLengthTheorem} (1) Let $P$ and $Q$ be vertex-colored posets with vertices colored by a set $I$.  Suppose that for each $i \in I$, the vertices of color $i$ in $Q$ coincide with the vertices of color $i$ in $P$.  Further suppose that $Q$ is a weak subposet of $P$.  Let $K := \mathbf{J}_{color}(P)$ and $L := \mathbf{J}_{color}(Q)$.  Then $K$ is a full-length edge-colored sublattice of $L$.\\ 
(2) Conversely, suppose $L$ is a diamond-colored distributive lattice with edges colored by a set $I$, with full-length sublattice $K$. Let $Q := \mathbf{j}_{color}(L)$ and $P := \mathbf{j}_{color}(K)$. 
\begin{itemize}
	\item For any join irreducible $\xelt$ in $L$, the set $\{\yelt \in K\, \mid \, x \leq_{L} \yelt\}$ has a unique minimal element $\welt_{\xelt}$, and $\welt_{\xelt}$ is a join irreducible in $K$. 
	\item The function $\phi:Q \longrightarrow P$ given by $\phi(\xelt) := \welt_{\xelt}$ is a vertex-color-preserving bijection, and if $\uelt \leq_{Q} \velt$ then $\phi(\uelt) \leq_{P} \phi(\velt)$. 
Now let $Q'$ be the set $P$ and declare that $\phi(\uelt) \leq_{Q'} \phi(\velt)$ if and only if $\uelt \leq_{Q} \velt$.  Then $Q'$ is a weak subposet of $P$ and $Q' \cong Q$ as vertex-colored posets.
\end{itemize}
\end{Thm}

\section{A classical example}\label{s:classicalL}
In this section we begin to explore a classical family of diamond-colored distributive lattices -- the aforementioned type $\myA$ fundamental lattices -- as an illustration of the methodology of the preceding section. 

\begin{figure}[H]
\begin{center}
\scalebox{.6}{


\begin{tikzpicture}

	\node at (-.75,0){$N-k$};
	\node at (-1.05,1){$N-k-1$};
	\node at (-.35,2){3};
	\node at (-.35,3){2};
	\node at (-.35,4){1};

	\node at (1.45,0){$N-k+1$};
	\node at (1.75,1){$N-k$};
	\node at (2.15,2){4};
	\node at (2.15,3){3};
	\node at (2.15,4){2};

	\node at (5.75,0){$N-1$};
	\node at (5.75,1){$N-2$};
	\node at (5.75,2){$k+2$};
	\node at (5.75,3){$k+1$};
	\node at (6.15,4){$k$};

	\foreach \i in {0,2.5,6.5}
	{
		\draw[dashed] (\i,1) -- (\i,2);
		\foreach \j in {0,2,3}
		{
				\node[circle, fill=black] at (\i,\j){};
				\draw (\i,\j) -- (\i,\j+1);
		}
		\foreach \j in {1,4}
			\node[circle, fill=black] at (\i,\j){};
	}
	
	\node[font=\LARGE] at (1.25,2){$\oplus$};
	\node[font=\LARGE] at (4,2){$\oplus \cdots \oplus$};
	
	\node[font=\Large] at (-3.5,2){$\widetilde{P}(k,N-k) :=$};

\end{tikzpicture} }
\end{center}
\caption{$\widetilde{P}(k,N-k)$.  The vertex colors are indicated.}
\label{fig:ChainSumFigure}
\end{figure}
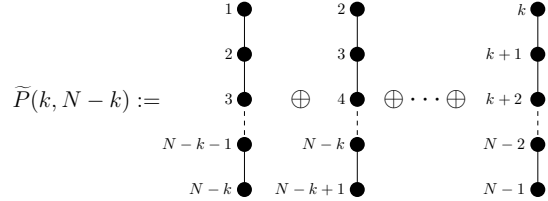

Throughout this section, $N$ is a fixed positive integer and $k$ is an integer with $1\leq k\leq N-1$. 
The poset $\widetilde{P}(k,N-k)$ is the disjoint sum of vertex-colored chains each with $N-k$ elements, as depicted in Figure \ref{fig:ChainSumFigure}. 

The chains are numbered from $1$ to $k$ reading from left to right.  In the $c$th chain, we number the vertices from $1$ to $N-k$ reading from bottom to top and color them from $c$ to $N-k-c+1$.  We may refer to a vertex by its \emph{row-column coordinates} $(r,c)$ for $1\leq r\leq k$ and $1\leq c\leq N$, yielding $VC(r,c) = N-k+r-c$. 

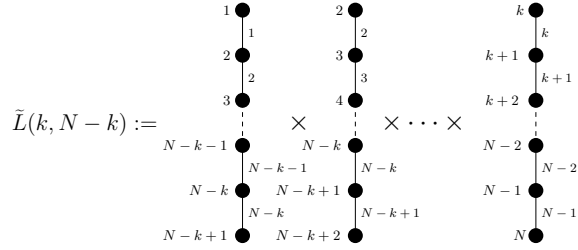
\begin{figure}[h]
\begin{center}
\scalebox{.6}{
%

\begin{tikzpicture}

	\node[font=\small] at (-1.05,0){$N-k+1$};
	\node[font=\small] at (-.75,1){$N-k$};
	\node[font=\small] at (-1.05,2){$N-k-1$};
	\node[font=\small] at (-.35,3){3};
	\node[font=\small] at (-.35,4){2};
	\node[font=\small] at (-.35,5){1};

	\node[font=\small] at (1.45,0){$N-k+2$};
	\node[font=\small] at (1.45,1){$N-k+1$};
	\node[font=\small] at (1.75,2){$N-k$};
	\node[font=\small] at (2.15,3){4};
	\node[font=\small] at (2.15,4){3};
	\node[font=\small] at (2.15,5){2};

	\node[font=\small] at (6.15,0){$N$};
	\node[font=\small] at (5.75,1){$N-1$};
	\node[font=\small] at (5.75,2){$N-2$};
	\node[font=\small] at (5.75,3){$k+2$};
	\node[font=\small] at (5.75,4){$k+1$};
	\node[font=\small] at (6.15,5){$k$};

	\foreach \i in {0,2.5,6.5}
	{
		\draw[dashed] (\i,2) -- (\i,3);
		\foreach \j in {0,...,5}
				\node[circle, fill=black] at (\i,\j){};
	}
	
	\node[font=\LARGE] at (1.25,2.5){$\times$};
	\node[font=\LARGE] at (4,2.5){$\times \cdots \times$};
	
	\node[font=\Large] at (-3.5,2.5){$\widetilde{L}(k,N-k) :=$};
	
	\draw (0,0) -- node[right, font=\footnotesize]{$N-k$} (0,1);
	\draw (0,1) -- node[right, font=\footnotesize]{$N-k-1$} (0,2);
	\draw (0,3) -- node[right, font=\footnotesize]{2} (0,4);
	\draw (0,4) -- node[right, font=\footnotesize]{1} (0,5);
	
	\draw (2.5,0) -- node[right, font=\footnotesize]{$N-k+1$} (2.5,1);
	\draw (2.5,1) -- node[right, font=\footnotesize]{$N-k$} (2.5,2);
	\draw (2.5,3) -- node[right, font=\footnotesize]{3} (2.5,4);
	\draw (2.5,4) -- node[right, font=\footnotesize]{2} (2.5,5);
	
	\draw (6.5,0) -- node[right, font=\footnotesize]{$N-1$} (6.5,1);
	\draw (6.5,1) -- node[right, font=\footnotesize]{$N-2$} (6.5,2);
	\draw (6.5,3) -- node[right, font=\footnotesize]{$k+1$} (6.5,4);
	\draw (6.5,4) -- node[right, font=\footnotesize]{$k$} (6.5,5);

\end{tikzpicture} }
\end{center}
\caption{The edge-colored distributive lattice $\widetilde{L}(k,N-k)$.  Vertex colors indicated on the left of each chain, edge colors on the right.}
\label{fig:ChainProduct}
\end{figure}

Declare $\widetilde{L}(k,N-k) := \mathbf{J}_{color}(\widetilde{P}(k,N-k))$.  We can view $\widetilde{L}(k,N-k)$ as the product of edge-colored chains as depicted in Figure \ref{fig:ChainProduct}.

So we can view $\widetilde{L}(k,N-k)$ as a collection of $k$-tuples, where the $i$th coordinate of each $k$-tuple is a choice of a label from the $i$th chain in the product of chains.  Writing each $k$-tuple as a column vector, we discern from Proposition \ref{p:FullLengthWithinProduct} (\ref{p:FullLengthWithinProduct1}) that 
\[\widetilde{L}(k,N-k) \cong \left\{\left(\begin{array}{c}{T}_{1}\\ {T}_{2}\\ \vdots\\ {T}_{k}\end{array}\right)\, \middle| \, \begin{array}{c}\mbox{For all $j$, ${T}_{j}\in\mathbb{Z}$, and}\\ \ j \leq {T}_{j} \leq N-k-1+j\end{array}\right\},\]
where the partial order of these column vectors is by ``reverse component-wise comparison,'' i.e. for such $k$-tuples ${S}$ and ${T}$ we have ${S} \leq {T}$ if and only if each ${S}_{j} \geq {T}_{j}$.  We have ${S} \xrightarrow{i} {T}$ if and only if each ${S}_{p} = {T}_{p}$ for $1 \leq p \leq k$ when $p \not= q$ while ${S}_{q} = i+1 = {T}_{q}+1$.  

It also follows that meets and joins of the column vectors are computed reverse component-wise as follows:
\[S \vee T = (S_{1} \wedge T_{1}, \ldots , S_{k} \wedge T_{k}) \mbox{\hspace*{0.25in} and \hspace*{0.25in}} S \wedge T = (S_{1} \vee T_{1}, \ldots , S_{k} \vee T_{k}).\] 
Moreover, \[\max(\widetilde{L}(k,N-k)) = \left(\begin{array}{c}1\\ 2\\ \vdots\\ k\end{array}\right) \mbox{\hspace*{0.25in} and \hspace*{0.25in}} \min(\widetilde{L}(k,N-k)) = \left(\begin{array}{c}N-k+1\\ N-k+2\\ \vdots\\ N\end{array}\right).\]

\paragraph{The ``column-strict'' tableau sublattice.} Now define a subset $\Lt(k,N-k)$ of $\widetilde{L}(k,N-k)$ as follows: 
\[\Ltab := \left\{\left(\begin{array}{c}{T}_{1}\\ {T}_{2}\\ \vdots\\ {T}_{k}\end{array}\right) \in \widetilde{L}(k,N-k)\ \, \middle| \, \begin{array}{c} 1 \leq {T}_{1} < {T}_{2} < \cdots < {T}_{k} \leq N\end{array}\right\},\]
and supply $\Ltab$ with the induced order from $\widetilde{L}(k,N-k)$. 
It is easy to check that Proposition \ref{p:FullLengthWithinProduct} (\ref{p:FullLengthWithinProduct2}) applies, so that $\Ltab$ is a full-length diamond-colored distributive sublattice of $\widetilde{L}(k,N-k)$. 

The superscript ``{\em tab}'' in the notation ``$\Ltab$'' indicates that the affiliated column vectors with strictly increasing entries (reading from top to bottom) can be viewed as column-shaped tableau. Tableaux often arise as objects that index bases for representations of certain algebraic structures, for example irreducible representations of the symmetric group or of the classical simple Lie algebras. The latter is the case here: the column vectors of $\Ltab$ index a natural weight basis for a certain fundamental representation of the type $\myA$ simple Lie algebra $sl(N,\mathbb{C})$.  In particular, if $\{e_{1}, e_{2}, \ldots , e_{N}\}$ is the standard basis for the defining complex representation $V$ of $sl(N,\mathbb{C})$, then $\{e_{T}\, |\, T \in  \Ltab\}$ is a weight basis for the $k$th exterior power $\bigwedge^{k}V$, where $e_{T}$ is the wedge product $e_{T_1} \wedge e_{T_2} \wedge \cdots \wedge e_{T_k}$ if $T$ is the column vector from $\Ltab$ with entries $T_{1}, T_{2}, \ldots , T_{k}$.  The so-called maximal vector of the representing space $\bigwedge^{k}V$ can be identified with the wedge product $e_{1} \wedge e_{2} \wedge \cdots \wedge e_{k}$ whose weight is the fundamental weight $\omega_{k}$. It can be seen that the edge-colored directed graph $\Ltab$ is the supporting graph for this ``wedge basis''  for the $k$th fundamental representation of $sl(N,\mathbb{C})$.  
It follows that a natural weight-generating function for $\Ltab$ coincides with the type $\myA$ Weyl bialternant $\chi_{_{\omega_{k}}}$ associated with the fundamental weight $\omega_{k}$; therefore $\Ltab$ is a splitting poset for the fundamental Weyl symmetric function $\chi_{_{\omega_{k}}}$. 
For further details, see \cite{Proc1}, \cite{proc3}, \cite{Don3}, \cite{HL}, and \cite{dds2}. 

In view of the fact that $\Ltab$ is a full-length diamond-colored distributive sublattice of $\widetilde{L}(k,N-k)$, we utilize Theorem \ref{t:FullLengthTheorem} to describe the vertex-colored poset of join irreducibles $P_{\myA}(k,N-k) := \mathbf{j}_{color}(\Ltab)$. 
The following shows how the row-column coordinates for an arbitrary vertex of $\widetilde{P}(k,N-k)$ corresponds with a join irreducible column $\xelt$ in $\widetilde{L}(k,N-k)$:
\[(r,c) \longleftrightarrow \left(\begin{array}{c}N-k+1\\ N-k+2\\ \vdots\\ N-k+r-c-1\\ N-k+r-c\\ N-k+r-c+1\\ \vdots \\ N\end{array}\right) =: \xelt,\]
which, in the notation of Theorem \ref{t:FullLengthTheorem}, corresponds to the following column in $\Ltab$:
\[\left(\begin{array}{c}N-k+1-c\\ N-k+2-c\\ \vdots\\ N-k+r-1-c\\ N-k+r-c\\ N-k+r+1-c\\ \vdots \\ N\end{array}\right)=: \welt_{\xelt}.\] 

Continuing in the notation of Theorem \ref{t:FullLengthTheorem}, we have $\phi(r,c) \leq \phi(s,d)$ in $P_{\myA}(k,N-k)$ if and only if $r \leq s$ and $c \leq d$.  
We refer to vertices of $P_{\myA}(k,N-k)$ using the row-column, or $(r,c)$ coordinates. It follows that the vertex-colored order diagram for $P_{\myA}(k,N-k)$ has the same edges as in $\widetilde{P}(k,N-k)$ plus some additional edges.  All this is to say that we have the following definition.

\begin{defn}\label{d:lA} Let $k$ and $N$ be integers with $1\leq k\leq N-1$.  Define the $k$th \emph{type} $\myA$ \emph{fundamental lattice}, $\La$, to be
\[\La := \mathbf{J}_{color}(P_{\myA}(k,N-k))\]
\end{defn}

\begin{figure}[H]
\begin{center}
\scalebox{.8}{
%

\begin{tikzpicture}

	\node[circle,fill=gray, draw=black] (a) at(1.5,-.5){};
	\node[circle,fill=gray, draw=black] (b) at(1,0){};
	\node[circle,fill=gray, draw=black] (c) at(2,0){};
	\node[circle,fill=gray, draw=black] (d) at(.5,.5){};
	\node[circle,fill=gray, draw=black] (e) at(1.5,.5){};
	\node[circle, draw=black] (f) at(2.5,.5){};
	\node[circle,fill=gray, draw=black] (g) at(0,1){};
	\node[circle,fill=gray, draw=black] (h) at(1,1){};
	\node[circle, draw=black] (i) at(2,1){};
	\node[circle, draw=black] (j) at(.5,1.5){};
	\node[circle, draw=black] (k) at(1.5,1.5){};
	\node[circle, draw=black] (l) at(1,2){};

	\draw (a) -- (b) -- (d) -- (g);
	\draw (c) -- (e) -- (h) -- (j);
	\draw (f) -- (i) -- (k) -- (l);
	\draw (a) -- (c) -- (f);
	\draw (b) -- (e) -- (i);
	\draw (d) -- (h) -- (k);
	\draw (g) -- (j) -- (l);

	\node at (3.5,.75){$\xrightarrow[\text{clockwise}]{135^\circ}$};

	\filldraw[fill=gray] (4.5,.5) rectangle (6,1.5);
	\filldraw[fill=gray] (6,1) rectangle (6.5,1.5);
	\foreach \j in {0,...,3}
		\draw (4.5,\j/2) -- (6.5, \j/2);
	\foreach \i in {9,...,13}
		\draw (\i/2,0) -- (\i/2,1.5);

\end{tikzpicture}
%
%
}
\end{center}
\caption{An order ideal from $P_{\mathsf{A}}(3,4)$ and its corresponding partition.}
\label{fig:rotation}
\end{figure}
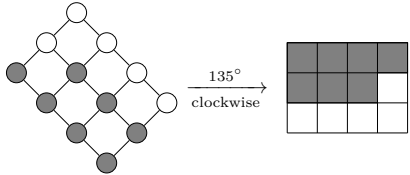

Note that $\La$ is a diamond-colored distributive lattice, colored by the set $[N-1]$. The elements of $\La$ are order ideals from $P_{\myA}(k,N-k)$, and computations within the environment of $\La$ are governed by Proposition \ref{p:EncapsulateProposition}. To make such computations concrete, it will help to identify an order ideal $\xelt$ from $P_{\myA}(k,N-k)$ with the sequence $(c_{1}(\xelt), c_{2}(\xelt), \ldots , c_{k}(\xelt))$, where for any $1 \leq r \leq k$ we have 
\[c_{r}(\xelt) := \left|\bigg\{c \in [N-k]\, \middle|\, (r,c) \in \xelt\bigg\}\right|.\] 
This order ideal can clearly be seen to be a partition, as in Figure \ref{fig:rotation}, since the partition corresponding to a given order ideal can be obtained by simply rotating a given order ideal clockwise $135^\circ$.  With this viewpoint, it is easy to see that the minimum element of $\La$ is the empty partition, and the maximum element is the full partition.  Furthermore, given two order ideals $\xelt$ and $\yelt$ viewed as partitions, we have that when one is overlaid on the other $\xelt\wedge\yelt$ is the partition consists of all boxes shaded in both $\xelt$ and $\yelt$ while $\xelt\vee\yelt$ contains any shaded box from $\xelt$ or $\yelt$.

\paragraph{Coordinatizations of $\L$.}  We now describe our four different coordinatizations of $\La$, or simply $\L$ when the parameters $k$ and $N$ are clear: partitions as $\Lp$ (as seen above, but denoted here for emphasis), tableaux as $\Lt$, circle diagrams as $\Lc$, and finally diagonal coordinates as $\Ld$.  Briefly, here are some of the advantages afforded by these different coordinatizations: We have seen that tableaux connect $\La$ to Lie algebra representation theory.  The type $\mathsf{A}$ fundamental lattices seem to have their origins as lattices of partitions, and this partition viewpoint helps us connect the type $\mathsf{A}$ fundamental lattices and the Domino Game.  Circle diagrams are useful for understanding the structure of ``two-color components'' of $\La$, and will lead to our most immediate connection between $\La$ and the Domino Game on $k\times(N-k)$ partitions.  The diagonal coordinates we present here seem to be new and afford a different, linear algebraic connection between $\La$ and the Domino Game.  

Note that because of the Domino Game and the simple connection to order ideals shown in Figure \ref{fig:rotation}, partitions will be our default mode of viewing the elements of $\L$.  Moreover, the other coordinatizations are best understood through this partition perspective.  In each case, the coordinatizations are edge-colored directed graphs with edges and edge colors induced by the universal object $\L$.   
 
\begin{defn}\label{def:lpart} Define $\mathscr{P}: \La \longrightarrow \{0,1,\ldots,N-k\}^{k}$ via $\mathscr{P}(\xelt) := (p_{r}(\xelt))_{r=1}^{k}$ where $p_{r}(\xelt) = c_{r}(\xelt)$.  Denote the image of $\mathscr{P}$ by $\Lpart$, or simply $\Lp$ when the dimensions are clear.  Here is a summary of the key data structures for $\Lp$: 
\begin{itemize}
	\item \textsc{Elements:} $\Lp\subseteq \BZ^k$, and $\sigma=(\sigma_1,\sigma_2,\ldots,\sigma_k)\in\BZ^k$ is in $\Lp$ $\Longleftrightarrow$ $N-k\geq \sigma_1\geq\sigma_2\geq\ldots\geq\sigma_k\geq 0$.
	\item \textsc{Directed edges:} For $\sigma,\tau\in\Lp$, we have a directed edge $\sigma\rightarrow\tau$ iff $\tau-\sigma=\varepsilon_l$, for some $l$ with $1\leq l\leq k$.
	\item \textsc{Edge colors:} Suppose $\sigma\rightarrow\tau$ in $\Lp$, so $\tau-\sigma=\varepsilon_l$ for some $l$, with $1\leq l\leq k$.  We say the edge is colored by $i$, and write $\sigma\xrightarrow{i}\tau$ iff $i=N-k-\tau_l+l$. It is clear that $i\in[N-1]$.
	\item \textsc{Visualization:} A $k \times (N-k)$ partition where the $l$th row contains $\gs_l$ left-justified shaded boxes. An edge of color $i$ exists between two partitions $\gs$ and $\gt$ if and only if one additional box is shaded when moving from the first diagram $\gs$ to the next, $\gt$ where $i=N-k-\tau_l+l$.
	
\begin{figure}[H]
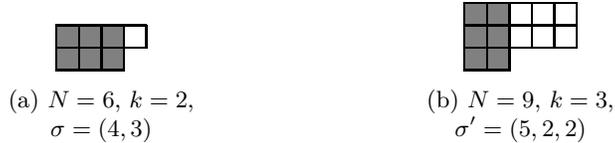

\centering
\begin{subfigure}[b]{0.45\textwidth}
\begin{center}
\ydiagram[*(white)]{0,3+1}*[*(gray)]{4,3}
\caption{$N=6$, $k=2$,\\ $\gs = (4,3)$}
\end{center}
\end{subfigure}
\begin{subfigure}[b]{0.45\textwidth}
\begin{center}
\ydiagram[*(white)]{0,2+3,2+3}*[*(gray)]{5,2,2}
\caption{$N=9$, $k=3$, \\ $\gs' = (5,2,2)$}
\end{center}
\end{subfigure}
\caption{Two partitions and their coordinatizations}
\end{figure}
	
\end{itemize}
\end{defn}

\begin{defn}\label{def:ltab} Define $\mathscr{T}: \La \longrightarrow [N]^{k}$ by the rule $\mathscr{T}(\xelt) := (T_{r}(\xelt))_{r=1}^{k}$ where $T_{r}(\xelt) := N-k+r-1-c_{r}(\xelt).$  Denote the image of $\mathscr{T}$ by $\Ltab$, or simply $\Lt$.  Here are the key data structures for $\Lt$:
\begin{itemize}
	\item \textsc{Elements:} $\Lt \subset \mathcal{P}([N])$, and $S \in \mathcal{P}([N])$ are in $\Lt$ $\Longleftrightarrow$ $|S|=k$. For $S\in\Lt$, it will be our convention to list the elements of the set in increasing order.  That is, for $S=\{S_1,S_2,\ldots,S_k\}$ we have $S_i<S_{i+1}$ for $i=1,\ldots, k-1$.
	\item \textsc{Edges:} For $S, T \in \Lt$, we have $S \rightarrow T$ if and only if $S = (T-\{i\})\cup \{i+1\}$ for some $i \in [N-1]$.
	\item \textsc{Edge colors:} For distinct subsets $S$ and $T$, we say $S\xrightarrow{i}T$ iff $S=(T-\{i\})\cup \{i+1\}$.
	\item \textsc{Visualization:} A $k \times 1$ tableau listing, from least to greatest, the number of unmarked boxes plus the row number, for each row in the partition. Alternately, we may attach an empty descending staircase to the right of the associated partition and count the number of unmarked boxes to obtain the same result. This method can be seen in Figure \ref{fig:LTabDiagrams} below. An edge of color $i$ exists between two tabular coordinatizations iff the value $i$ increases to $i+1$ from one diagram to the next.

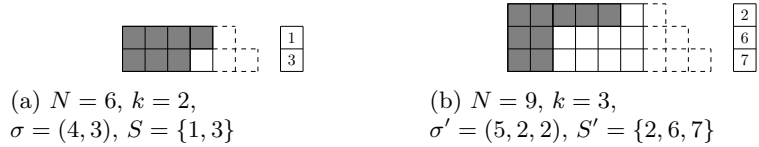
\begin{figure}[H]
\centering
\begin{subfigure}[b]{0.45\textwidth}
\begin{center}
\scalebox{.6}{
\begin{tikzpicture}
	\filldraw[fill=gray] (0,0) rectangle (1.5,1);
	\filldraw[fill=gray] (1.5,.5) rectangle (2,1);

	\foreach \i in {0,...,2}
		\draw (0,\i/2) -- (2,\i/2);
	\foreach \j in {0,...,4}
		\draw (\j/2,0) -- (\j/2,1);
	\draw[dashed] (2,0) -- (3,0);
	\draw[dashed] (2,.5) -- (3,.5);
	\draw[dashed] (2,1) -- (2.5,1);
	\draw[dashed] (2.5,0) -- (2.5,1);
	\draw[dashed] (3,0) -- (3,.5);
	
	\draw (3.5,0) -- (4,0);
	\draw (3.5,.5) -- (4,.5);
	\draw (3.5,1) -- (4,1);
	\draw (3.5,0) -- (3.5,1);
	\draw (4,0) --(4,1);
	
	\node at (3.75,.25){3};
	\node at (3.75,.75){1};
	
\end{tikzpicture}
}
\caption{$N=6$, $k=2$,\\ $\gs = (4,3)$, $S=\{1,3\}$}
\end{center}
\end{subfigure}
\begin{subfigure}[b]{0.45\textwidth}
\begin{center}
\scalebox{.6}{
\begin{tikzpicture}
	\filldraw[fill=gray] (0,0) rectangle (1,1.5);
	\filldraw[fill=gray] (1,1) rectangle (2.5,1.5);

	\foreach \i in {0,...,3}
		\draw (0,\i/2) -- (3,\i/2);
	\foreach \j in {0,...,6}
		\draw (\j/2,0) -- (\j/2,1.5);
	\draw[dashed] (3,0) -- (4.5,0);
	\draw[dashed] (3,.5) -- (4.5,.5);
	\draw[dashed] (3,1) -- (4,1);
	\draw[dashed] (3,1.5) -- (3.5,1.5);
	\draw[dashed] (3.5,0) -- (3.5,1.5);
	\draw[dashed] (4,0) -- (4,1);
	\draw[dashed] (4.5,0) -- (4.5,.5);
	
	\draw (5,0) -- (5.5,0);
	\draw (5,.5) -- (5.5,.5);
	\draw (5,1) -- (5.5,1);
	\draw (5,1.5) -- (5.5,1.5);
	\draw (5,0) -- (5,1.5);
	\draw (5.5,0) --(5.5,1.5);
	
	\node at (5.25,.25){7};
	\node at (5.25,.75){6};
	\node at (5.25,1.25){2};
	
\end{tikzpicture}
}
\caption{$N=9$, $k=3$, \\ $\gs' = (5,2,2)$, $S'=\{2,6,7\}$}
\end{center}
\end{subfigure}
\caption{\label{fig:LTabDiagrams}%
Two partitions with augmented staircases and their associated tableaux.}
\end{figure}
\end{itemize}
\end{defn}

\begin{defn}\label{def:lcirc} Define $\mathscr{B}: \La \longrightarrow \{0,1\}^{N}$ via $\mathscr{B}(\xelt) := (b_{i}(\xelt))_{i=1}^{N}$ where 
\[b_{i}(\xelt) := \left\{\begin{array}{cl}1 & \hspace*{0.0in}\mbox{if $i$ is an entry in $\mathscr{T}(\xelt)$}\\ 0 & \hspace*{0.0in}\mbox{otherwise}\end{array}\right.\] 
Denote the image of $\mathscr{B}$ by $\Lcirc$, or simply $\Lc$.  Here are the key data structures for $\Lc$:
\begin{itemize}
	\item \textsc{Elements:} $\Lc\subseteq\BZ_2^N$, and $s=(s_1,s_2,\ldots,s_N)\in\BZ_2^{N}$ $\Longleftrightarrow$ $\sum s_i=k$.
	\item \textsc{Edges:} For $s,t\in\Lc$, we have $s\rightarrow t$ iff $t-s=-\varepsilon_l+\varepsilon_{l+1}$ for some $l$ with $1\leq l\leq N-1$, where $\varepsilon_l$ is a binary sequence of appropriate length with a 1 in position $l$ and 0's elsewhere.  
	\item \textsc{Edge colors:} Suppose $s\rightarrow t$ in $\Lc$, so $t-s=-\varepsilon_l+\varepsilon_{l+1}$ for some $l$ with $1\leq l\leq N-1$.  Write $s\xrightarrow{i}t$ iff $i=l$.  Note: It is clear that $i\in[N-1]$.
	\item \textsc{Visualization:} A $2 \times \left\lceil N/2 \right\rceil$ tableau (with the two rightmost boxes joined into one tall box when $N$ is odd) numbered as in Figure \ref{fig:LCircleDiagrams}. For $N$ even, we number the boxes clockwise, beginning in the top left box. For $N$ odd, we number the boxes counter-clockwise, beginning in the lower left corner.  The box labeled $i$ is marked with a dot iff $s_i=1$, $1 \leq i\leq N$. An edge of color $l$, denoted as $c_l$ below, exists between two circle diagrams if there is a single dot movement across the ``wall'' from box $l$ and to box $l+1$.  

\begin{figure}[H]
\centering
\begin{subfigure}[b]{0.45\textwidth}
\begin{center}
\begin{tikzpicture}

\foreach \place/\name in {{(0,.5)/ta}, {(.5,.5)/tb}, {(1,.5)/tc}, {(1.5,.5)/td}, {(2,.5)/te}, {(2.5,.5)/tf}, {(0,-.5)/ba}, {(.5,-.5)/bb}, {(1,-.5)/bc}, {(1.5,-.5)/bd}, {(2,-.5)/be}, {(2.5,-.5)/bf}, {(0,0)/ma},  {(1.5,0)/md}, {(2,0)/me}, {(2.5,0)/mf}}
    \node[circle, fill=white] (\name) at \place {};

\draw[line width=.8pt] (ta.base) -- (td.base);
\draw[line width=.8pt] (te.base) -- (tf.base);
\draw[line width=.8pt] (ma.base) -- (md.base);
\draw[line width=.8pt] (me.base) -- (mf.base);
\draw[line width=.8pt] (ba.base) -- (bd.base);
\draw[line width=.8pt] (be.base) -- (bf.base);
\draw[line width=.8pt] (ta.base) -- (ba.base);
\draw[line width=.8pt] (tb.base) -- (bb.base);
\draw[line width=.8pt] (tc.base) -- (bc.base);
\draw[line width=.8pt] (td.base) -- (bd.base);
\draw[line width=.8pt] (te.base) -- (be.base);
\draw[line width=.8pt] (tf.base) -- (bf.base);

\node [font=\scriptsize] at (0.25,.7){1};
\node [font=\scriptsize] at (.75,.7){2};
\node [font=\scriptsize] at (1.25,.7){3};
\node [font=\tiny] at (1.75,.7){$\cdots$};
\node [font=\scriptsize] at (2.25,.7){$\frac{N}{2}$};

\node [font=\tiny] at (0.22,-.73){$N$};
\node [font=\tiny] at (.75,-.75){$N\!\!-\!\!1$};
\node [font=\tiny] at (1.25,-.75){$N\!\!-\!\!2$};
\node [font=\tiny] at (1.75,-.75){$\cdots$};
\node [font=\scriptsize] at (2.25,-.75){$\frac{N}{2}\!\!+\!\!1$};

\node [font=\tiny] at (1.75,-0.25){$\cdots$};
\node [font=\tiny] at (1.75,0.25){$\cdots$};

\draw[red, dashed] (.5,.05) -- (.5,.85) node[above, font=\footnotesize]{$c_1$};
\draw[blue, dashed] (1,.05) -- (1,.85) node[above, font=\footnotesize]{$c_2$};
\draw[green, dashed] (1.5,.05) -- (1.5,.85) node[above, font=\footnotesize]{$c_3$};
\draw[gold,dashed] (2,.05) -- (2,.85) node[above, font=\footnotesize]{$c_\subsc{$\frac{N}{2}\!\!-\!\!1$}$};
\draw[seafoam, dashed] (2.05,0) -- (2.75,0) node[right, font=\footnotesize]{$c_\subsc{$\frac{N}{2}$}$};
\draw[violet, dashed] (2,-.05) -- (2,-.85) node[below, font=\footnotesize]{$\;\;c_\subsc{$\!\frac{N}{2}\!\!+\!\!1$}$};
\draw[cyan, dashed] (1.5,-.05) -- (1.5,-.85) node[below, font=\footnotesize]{$c_\subsc{N-3}$};
\draw[orange, dashed] (1,-.05) -- (1,-.85) node[below, font=\footnotesize]{$c_\subsc{N-2}$};
\draw[magenta, dashed] (.5,-.05) -- (.5,-.85) node[below, font=\footnotesize]{$c_\subsc{N-1}$};

\end{tikzpicture}%
\end{center}
\caption{Circle diagram for $N$ even}
\end{subfigure}
\begin{subfigure}[b]{0.45\textwidth}
\begin{center}
%



\begin{tikzpicture}

%
%

\foreach \i in {0,...,6}
	\draw[line width=.8pt] (\i/2,0) -- (\i/2,1);

\foreach \j in {0,...,2}
	\draw[line width=.8pt] (0,\j/2) -- (1.5,\j/2);
\draw[line width=.8pt] (2,0) -- (3,0);
\draw[line width=.8pt] (2,.5) -- (2.5,.5);
\draw[line width=.8pt] (2,1) -- (3,1);

\draw[red, dashed] (.5,.45) -- (.5,-.35) node[below, font=\footnotesize]{$c_1$};
\draw[blue, dashed] (1,.45) -- (1,-.35) node[below, font=\footnotesize]{$c_2$};
\draw[green, dashed] (1.5,.45) -- (1.5,-.35) node[below, font=\footnotesize]{$c_3$};
\draw[gold,dashed] (2,.45) -- (2,-.35) node[below, font=\footnotesize]{$\;\;\;c_\subsc{$\left\lfloor\!\!\frac{N}{2}\!\!\right\rfloor\!\!-\!\!1$}$};
\draw[seafoam, dashed] (2.5,.45) -- (2.5,0);
\draw[seafoam, dashed] (2.5,0) -- (2.75,-.35) node[below, font=\footnotesize]{$c_\subsc{$\left\lfloor\!\!\frac{N}{2}\!\!\right\rfloor$}$};
\draw[cyan, dashed] (2,.55) -- (2,1.35) node[above, font=\footnotesize]{$c_\subsc{$\left\lceil\!\!\frac{N}{2}\!\!\right\rceil\!\!+\!\!1$}$};
\draw[violet, dashed] (2.5,.55) -- (2.5,1);
\draw[violet, dashed] (2.5,1) -- (2.75,1.35) node[above, font=\footnotesize]{$c_\subsc{$\left\lceil\!\!\frac{N}{2}\!\!\right\rceil$}$};
\draw[magenta, dashed] (1.5,.55) -- (1.5,1.35) node[above, font=\footnotesize]{$c_\subsc{N-3}$};
\draw[orange, dashed] (1,.55) -- (1,1.35) node[above, font=\footnotesize]{$c_\subsc{N-2}$};
\draw[lavender, dashed] (.5,.55) -- (.5,1.35) node[above, font=\footnotesize]{$c_\subsc{N-1}$};

\node[font=\tiny] at (0.25,1.2){$N$};
\node[font=\tiny] at (.75,1.2){$N\!\!-\!\!1$};
\node[font=\tiny] at (1.25,1.2){$N\!\!-\!\!2$};
\node[font=\tiny] at (1.75,1.2){$\cdots$};
\node[font=\tiny] at (2.2,1.2){$\;\;\left\lceil\!\!\frac{N}{2}\!\!\right\rceil\!\!+\!\!1$};
\node[font=\tiny] at (3.25,.5){$\left\lceil\!\!\frac{N}{2}\!\!\right\rceil$};

\node[font=\scriptsize] at (0.25,-.25){$1$};
\node[font=\scriptsize] at (.75,-.25){$2$};
\node[font=\scriptsize] at (1.25,-.25){$3$};
\node[font=\tiny] at (1.75,-.25){$\cdots$};
\node[font=\tiny] at (2.25,-.25){$\left\lfloor\!\!\frac{N}{2}\!\!\right\rfloor$};

\node[font=\tiny] at (1.75,0.25){$\cdots$};
\node[font=\tiny] at (1.75,0.75){$\cdots$};

\end{tikzpicture}
%
\end{center}
\caption{Circle diagram for $N$ odd}
\end{subfigure}

\caption{\label{fig:LCircleDiagrams}%
Generic circle diagram labellings.}

\end{figure}
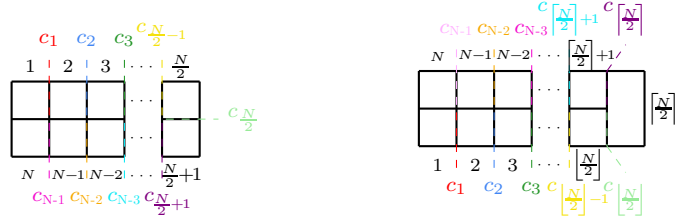

\begin{figure}[H]
\centering
\begin{subfigure}[b]{0.45\textwidth}
\begin{center}
\scalebox{.6}{%
\begin{tikzpicture}

\foreach \place/\name in {{(0,.5)/ta}, {(.5,.5)/tb}, {(1,.5)/tc}, {(1.5,.5)/td}, {(0,-.5)/ba}, {(.5,-.5)/bb}, {(1,-.5)/bc}, {(1.5,-.5)/bd}, {(0,0)/ma},  {(1.5,0)/md}, {(1.5,0)/md}}
    \node[circle, fill=white] (\name) at \place {};

\node[circle, fill=black] (dot1) at (.25,.25){};
\node[circle, fill=black] (dot2) at (1.25,-.25){};

\draw[line width=1pt] (ta.base) -- (td.base);
\draw[line width=1pt]  (ma.base) -- (md.base);
\draw[line width=1pt]  (ba.base) -- (bd.base);
\draw[line width=1pt]  (ta.base) -- (ba.base);
\draw[line width=1pt]  (tb.base) -- (bb.base);
\draw[line width=1pt]  (tc.base) -- (bc.base);
\draw[line width=1pt]  (td.base) -- (bd.base);

\draw[dashed,red] (.5,.05) -- (.5,.9)node[above, font=\footnotesize]{$c_1$};
\draw[dashed,blue] (1,.05) -- (1,.9)node[above, font=\footnotesize]{$c_2$};
\draw[dashed,green] (1.05,0) -- (1.8,0)node[right, font=\footnotesize]{$c_3$};
\draw[dashed,violet] (1,-.05) -- (1,-.9)node[below, font=\footnotesize]{$c_4$};
\draw[dashed,orange] (.5,-.05) -- (.5,-.9)node[below, font=\footnotesize]{$c_5$};

\node (v1) at (0.25,.7){1};
\node (v2) at (.75,.7){2};
\node (v3) at (1.25,.7){3};

\node (v6) at (0.25,-.7){6};
\node (v7) at (.75,-.7){5};
\node (v8) at (1.25,-.7){4};

\end{tikzpicture}%
}
\end{center}
\caption{$N=6$, $k=2$,\\ $\gs = (4,3)$, $S=\{1,3\}$}
\end{subfigure}
\begin{subfigure}[b]{0.45\textwidth}
\begin{center}
\scalebox{.6}{%

\begin{tikzpicture}

\foreach \place/\name in {{(0,.5)/ta}, {(.5,.5)/tb}, {(1,.5)/tc}, {(1.5,.5)/td}, {(2,.5)/te}, {(2.5,.5)/tf}, {(0,-.5)/ba}, {(.5,-.5)/bb}, {(1,-.5)/bc}, {(1.5,-.5)/bd}, {(2,-.5)/be}, {(2.5,-.5)/bf}, {(0,0)/ma}, {(2,0)/me}}
    \node[circle, fill=white] (\name) at \place {};

\node[fill=white] (tl) at (2.1,.5) {};
\node[fill=white] (tr) at (2.6,.5) {};
\node[fill=white] (bl) at (2.1,0) {};
\node[fill=white] (br) at (2.6,0) {};

\node[circle, fill=black] (dot1) at (.75,-.25) {};
\node[circle, fill=black] (dot2) at (1.25,.25){};
\node[circle, fill=black] (dot3) at (1.75,.25){};

\draw[line width=1pt] (ta.base) -- (tf.base);
\draw[line width=1pt] (ma.base) -- (me.base);
\draw[line width=1pt] (ba.base) -- (bf.base);
\draw[line width=1pt] (ta.base) -- (ba.base);
\draw[line width=1pt] (tb.base) -- (bb.base);
\draw[line width=1pt] (tc.base) -- (bc.base);
\draw[line width=1pt] (td.base) -- (bd.base);
\draw[line width=1pt] (te.base) -- (be.base);
\draw[line width=1pt] (tf.base) -- (bf.base);

\draw[dashed,red] (.5,-.05) -- (.5,-.9)node[below, font=\footnotesize]{$c_1$};
\draw[dashed,blue] (1,-.05) -- (1,-.9)node[below, font=\footnotesize]{$c_2$};
\draw[dashed,green] (1.5,-.05) -- (1.5,-.9)node[below, font=\footnotesize]{$c_3$};
\draw[dashed,violet] (2,-.05) -- (2,-.9)node[below, font=\footnotesize]{$c_4$};
\draw[dashed,orange] (2,.05) -- (2,.9)node[above, font=\footnotesize]{$c_5$};
\draw[dashed,cyan] (1.5,.05) -- (1.5,.9)node[above, font=\footnotesize]{$c_6$};
\draw[dashed,gold] (1,.05) -- (1,.9)node[above, font=\footnotesize]{$c_7$};
\draw[dashed,magenta] (.5,.05) -- (.5,.9)node[above, font=\footnotesize]{$c_8$};

\node (v1) at (0.25,.7){9};
\node (v2) at (.75,.7){8};
\node (v3) at (1.25,.7){7};
\node (v4) at (1.75,.7){6};
\node (v5) at (2.7,0){5};

\node (v6) at (0.25,-.7){1};
\node (v7) at (.75,-.7){2};
\node (v8) at (1.25,-.7){3};
\node (v9) at (1.75,-.7){4};

\end{tikzpicture}

}
\end{center}
\caption{$N=9$, $k=3$, \\ $\gs' = (5,2,2)$, $S'=\{2,6,7\}$}
\end{subfigure}

\caption{
Two circle diagrams for the given tabular coordinatizations.}

\end{figure}
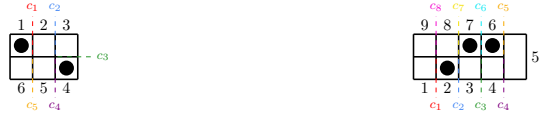

\end{itemize}
\end{defn} 

\begin{figure}[H]
\centering
\begin{subfigure}[b]{0.45\textwidth}
\begin{center}
\scalebox{.6}{%
\begin{tikzpicture}
	\filldraw[fill=gray] (0,0) rectangle (1.5,1);
	\filldraw[fill=gray] (1.5,.5) rectangle (2,1);

	\foreach \i in {0,...,2}
		\draw (0,\i/2) -- (2,\i/2);
	\foreach \j in {0,...,4}
		\draw (\j/2,0) -- (\j/2,1);
	\draw[dashed] (1.25,1.25) -- (2,.5);
	\foreach \i in {2,...,4}
		\draw[dashed] (\i/2-1.25 ,1.25) -- (\i/2,0);
	\draw[dashed] (-.25,.75) -- (.5,0);
	
\end{tikzpicture}%
}
\end{center}
\caption{$N=6$, $k=2$,\\ $\gs = (4,3)$, $\ds=(1,1,2,2,1)$}
\end{subfigure}
\begin{subfigure}[b]{0.45\textwidth}
\begin{center}
\scalebox{.6}{%
\begin{tikzpicture}
	\filldraw[fill=gray] (0,0) rectangle (1,1.5);
	\filldraw[fill=gray] (1,1) rectangle (2.5,1.5);

	\foreach \i in {0,...,3}
		\draw (0,\i/2) -- (3,\i/2);
	\foreach \j in {0,...,6}
		\draw (\j/2,0) -- (\j/2,1.5);
	\draw[dashed] (1.75,1.75) -- (3,.5);
	\draw[dashed] (2.25,1.75) -- (3,1);
	\foreach \i in {3,...,6}
		\draw[dashed] (\i/2-1.75 ,1.75) -- (\i/2,0);
	\draw[dashed] (-.25,.75) -- (.5,0);
	\draw[dashed] (-.25,1.25) -- (1,0);
	
\end{tikzpicture}
}
\end{center}
\caption{$N=9$, $k=3$, \\ $\gs' = (5,2,2)$, $\ds'=(0,1,1,1,1,2,2,1)$}
\end{subfigure}

\caption{\label{fig:LDiagDiagrams}%
Diagonal examples with their associated partitions.}

\end{figure}
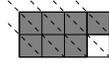
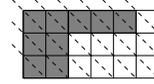

\begin{defn}\label{def:ldiag} Define $\mathscr{D}: \La \longrightarrow \{0,1,\ldots,k\}^{N-1}$ via $\mathscr{D}(\xelt) := (d_{i}(\xelt))_{i=1}^{N-1}$ where 
\[d_i(\xelt):=\begin{cases}
	\max\{0,r\mid c_r(\xelt)\geq (N-k)-(i-j)\}_{r=1}^{m_i} & \text{ if } 1\leq i\leq N-k\\
	\max\{0,r\mid c_{(r+i)-(N-k)}(\xelt)\geq r\}_{r=1}^{n_i} & \text{ if } N-k \leq i\leq N-1.
	\end{cases}\]
and where $m_i=\min\{i,k\}$, and if $N-k\leq i\leq N-1$, put $n_i=\min\{N-k,N-i\}$.  Denote the image of $\mathscr{D}$ by $\Ldiag$, or simply $\Ld$.  Note that this map is simply counting the shaded dots in vertical alignment from left to right in the order ideal $\xelt$ when pictured as in Figure \ref{fig:rotation} (cf. Equation \ref{e:gampd} in Definition \ref{d:part-diag}).  Here is a summary of the key data structures for $\Ld$: 

\begin{itemize}
	\item \textsc{Elements:} $\Ld\subseteq\BZ^{N-1}$, and $\ds=(\ds_1,\ds_2,\ldots,\ds_{N-1})\in\Ld$ iff
	\begin{itemize}
		\item $\ds_i\leq \min\{i,k\}$, for $1\leq i\leq N-k$,
		\item $\ds_i\leq \min\{N-i, N-k\}$, for $N-k\leq i\leq N-1$, 
		\item For $1\leq i\leq N-k-1$, $\ds_{i+1}=\ds_i$ or $\ds_i+1$, and
		\item For $N-k\leq i\leq N-2$, $\ds_{i+1}=\ds_i$ or $\ds_i-1$.
	\end{itemize}
	\item \textsc{Edges:} For $\ds,\dt\in\Ld$, we have $\ds\rightarrow \dt$, iff $\dt-\ds=\varepsilon_l$ for some $l$ with $1\leq l\leq N-1$.
	\item \textsc{Edge colors:} Suppose $\ds\rightarrow \dt$ in $\Ld$, so $\dt-\ds=\varepsilon_l$ for some $l$ with $1\leq l\leq N-1$.  Then $\ds\xrightarrow{i}\dt$ iff $i=l$.
	\item \textsc{Visualization:} The number of marked/shaded boxes of a given diagonal in the associated partition diagram where the diagonals slant from northwest to southeast, starting with the rightmost box in the top row of the partition. Note that the first $N-k$ diagonals are viewed across the top of the partition and the remaining $k-1$ move down the left side of the diagram. See Figure \ref{fig:LDiagDiagrams} for two examples.
\end{itemize}
\end{defn}

All of this is to say that $\Lt$, $\Lp$, $\Lc$ and $\Ld$ are truly coordinatizations of $\L$, in that they are isomorphic descriptions of the same combinatorial object: A diamond-colored distributive lattice colored by the set $[N-1]$.  For a complete example of these coordinatizations see Figure \ref{fig:L24Full}.

\section{Move-minimizing games}\label{s:games}
In this section we consider certain kinds of combinatorial one-player games and connect them to diamond-colored modular and distributive lattices. 

Suppose a single player is given a finite set of objects which are allowed to be manipulated by a certain finite set of moves, with objects and moves described combinatorially.  
Here are the ground rules: 
\begin{itemize}
	\item The player can ``legally'' apply a given move only to some known subset of the objects; 
	\item when the player legally applies the move to a given object $\xelt$, the outcome could be any one of a known subset of objects, and s/he chooses one of these resulting objects $\yelt$; 
	\item and a given move can be reversed in the sense that if $\yelt$ is obtained by legally applying a move to $\xelt$, then we also allow the player to obtain $\xelt$ from $\yelt$.
\end{itemize}
So, given two objects $\selt$ and $\telt$, the primary objective of the game is to minimize the number of moves needed to obtain $\telt$ from $\selt$.   Within this move-minimizing context, it will make sense to omit the following possibilities: when a move is legally applied to some object $\xelt$, then $\xelt$ cannot be among the outcomes, and when some object $\yelt$ is among the outcomes obtained by legally applying some move to $\xelt$, then $\yelt$ is not among the outcomes when any other move is legally applied to $\xelt$. 

For example, consider the Towers of Hanoi game with $n$ distinctly-sized disks arranged in stacks on three posts (say, post$\!$ \#1, post$\!$ \#2, and post$\!$ \#3), with the usual restriction that a larger disk cannot be stacked on top of a smaller disk.  
It is well-known that there are $3^{n}$ possible configurations of the disks; these configurations will be our objects.  
We could consider three possible moves, with each move corresponding to an ordered pair $(i,j)$ with $i, j \in \{1,2,3\}$ and $i < j$, defined as the disk at the top of the stack on post$\!$ \#$i$ is taken to the top of the stack on post$\!$ \#$j$.   
A goal of a move-minimizing game would be to obtain a given configuration from some other given configuration using the least number of moves.  
The undirected and uncolored graph underlying this move-minimizing game is well-known; analyses of shortest paths in this graph can be found in \cite{chan}, \cite{hinz}, and \cite{rom}. 

The alert reader notice that what we have described as a combinatorial one-player game of minimizing moves between objects coincides with analyzing shortest paths between vertices in some kind of edge-colored directed graph. 
We make this more precise below. 

\paragraph{Move-minimizing game digraphs.}\label{par:digraphs}  Let $\mathcal{G}$ be a finite, simple digraph with edges colored by some finite index set $I$. 
The vertices of $\mathcal{G}$ are to be regarded as objects while the directed and colored edges of $\mathcal{G}$ correspond to legal moves.   

The fundamental question for a \emph{move-minimizing game} is given any two objects $\selt$ and $\telt$ of such a graph $\mathcal{G}$, ``Is there a path in $\mathcal{G}$ from $\selt$ to $\telt$?''  If such a path exists, then we will consider three additional questions:  
\begin{enumerate}  
\item[(1)] What is an explicit formula for the distance from $\selt$ to $\telt$? 
\item[(2)] For any given color $i \in I$, what is an explicit formula for the minimum number of color $i$ edges amongst all shortest paths from $\selt$ to $\telt$? 
\item[(3)] For any given color $i \in I$, can a shortest path from $\selt$ to $\telt$ that minimizes the number of color $i$ edges be explicitly prescribed in terms of $\selt$ and $\telt$?
\end{enumerate} 

Call $\mathcal{G}$ the {\em move-minimizing game digraph} associated with this move-minimizing game. 
Notice that if we are not interested in using colors to distinguish between moves, then we may simply consider the uncolored directed graph, in which case question (2) coincides with question (1) while question (3) may be rephrased simply as ``Can a shortest path from $\selt$ to $\telt$ be explicitly prescribed in terms of the objects $\selt$ and $\telt$?'' Thus, answers to our four move-minimizing questions within the context of colored moves will also serve as answers to the ``uncolored'' versions of these questions. 

We understand that the meaning of the term ``explicit'' in questions (1), (2), and (3) is open to interpretation and certainly depends upon the combinatorial description of the objects in question. The following theorem connects move-minimizing games and diamond-colored modular and distributive lattices and indicates the level of explicitness we aim to achieve in answers to questions (1), (2), and (3). 

\begin{Thm}\label{t:MoveMinGameTheorem} Suppose a move-minimizing game digraph $\mathcal{G}$ is the order diagram for some diamond-colored modular lattice $L$ whose rank function is $\rho$.  
Let $\selt, \telt \in \mathcal{G}$.  If $L$ is distributive, then let $P := \mathbf{j}_{\mbox{\em color}}(L)$ be the vertex-colored poset of join irreducibles affiliated with $L$, and identify $\selt$ and $\telt$ as order ideals from $P$, so in particular $\selt$ and $\telt$ are subsets of $P$.  Then:\\
\begin{enumerate}
\setcounter{enumi}{-1}
	\item There exists a path in $\mathcal{G}$ from $\selt$ to $\telt$.
	\item The distance from $\selt$ to $\telt$ is $2\rho(\selt \vee \telt) - \rho(\selt) - \rho(\telt) = \rho(\selt) + \rho(\telt) - 2\rho(\selt \wedge \telt)$. If $L$ is distributive, then the quantities here may be determined as follows: $\rho(\selt) = |\selt|$, $\rho(\telt) = |\telt|$, $\rho(\selt \vee \telt) = |\selt \cup \telt|$, and $\rho(\selt \wedge \telt) = |\selt \cap \telt|$.
	\item For any given $i \in I$, any two shortest paths from $\selt$ to $\telt$ have the same number of color $i$ edges.  If $L$ is distributive, then this number is the number of color $i$ vertices in the set $(\selt \cup \telt) - \selt$ plus the number of color $i$ vertices in the set $(\selt \cup \telt) - \telt$, or equivalently the number of color $i$ vertices in $\selt - (\selt \cap \telt)$ plus the number of color $i$ vertices in $\telt - (\selt \cap \telt)$.
	\item Suppose now that $L$ is distributive.  Let $\ell_{1} := |\selt \cup \telt| - |\selt|$ and $\ell_{2} := |\selt \cup \telt| - |\telt|$.  Let $\xelt_{0} := \selt$.  For $1 \leq j \leq \ell_{1}$, let $u$ be any minimal element of the set $(\selt \cup \telt) - \xelt_{j-1}$ and declare $\xelt_{j} := \xelt_{j-1} \cup \{u\}$.  Let $\xelt_{\ell_{1}+\ell_{2}} := \telt$.  For $1 \leq j \leq \ell_{2}$, let $u$ be any minimal element of the set $(\selt \cup \telt) - \xelt_{\ell_{1} + j}$ and declare $\xelt_{\ell_{1} + j - 1} := \xelt_{\ell_{1} + j} \cup \{u\}$.  Then for $0 \leq j \leq \ell_{1}+\ell_{2}$, $\xelt_{j}$ is an order ideal from $P$, and $(\xelt_{0}, \xelt_{1}, \ldots , \xelt_{\ell_{1}+\ell_{2}})$ is a shortest path from $\selt$ to $\telt$ with $\xelt_{\ell_{1}} = \selt \cup \telt$. Similarly obtain a shortest path $(\selt = \yelt_{0}, \yelt_{1}, \ldots , \yelt_{k_{1}+k_{2}}=\telt)$ from $\selt$ to $\telt$ with $\yelt_{k_{1}} = \selt \cap \telt$, where $k_{1} = |\selt| - |\selt \cap \telt|$ and $k_{2} = |\telt| - |\selt \cap \telt|$.
\end{enumerate}
\end{Thm}
\begin{proof} The claims follow from Theorem \ref{t:ModularLatticeTheorem}, Lemma \ref{l:ColorsLemma}, and Proposition \ref{p:EncapsulateProposition}.
\end{proof} 

We close this section with some remarks. 
\begin{itemize}
\item The phrases ``God's algorithm'' and ``God's number'' are most often used in connection with a move-minimizing game modeled using the Cayley graph of some group related to the underlying puzzle, perhaps most famously in connection with Rubik's Cube (see for example \cite{joyn}).  In a similar sense, then, the preceding theorem determines God's number and God's algorithm for a move-minimizing game between two objects whenever the game digraph is the order diagram for a diamond-colored distributive lattice. 
\item At first glance it may seem unlikely that any interesting games could be modeled in this way. Here and in sequel papers we will showcase some move-minimizing games -- our so-called ``domino games'' -- that serve as nice illustrations of this approach. Furthermore, we will explore some of the many enumerative and algebraic contexts within which our move-minimizing domino game digraphs naturally occur. 
\end{itemize}

\section{The Domino Game digraph}\label{s:dominos}
We now fully describe the Domino Game including its description as a digraph with edges arising from legal domino moves.  Fix integers $N$ and $k$ with $1 \leq k \leq N-1$ and consider $k\times(N-k)$ partitions as the vertices of a graph we denote by $\Da$, or $\D$ when the parameters $k$ and $N$ are clear.  We color the game board checkerboard-style, with the upper right most square being red, and non-red squares being white.  The use of colors on the game board is a refinement of the game as described in the Introduction, and will yield a digraph that we can connect to $\La$  

Given any two $k\times(N-k)$ partitions (or ``legal shapes'') $\sigma$ and $\tau$, form a directed edge $\sigma\rightarrow\tau$ if $\tau$ can be obtained from $\sigma$ by any one of the following legal domino moves:
\begin{enumeratei}
	\item Remove from $\sigma$ (or add to $\tau$) \ydiagram[*(white)]{1+1}*[*(red!50)]{1}, \ydiagram[*(white)]{1,0}*[*(red!50)]{0,1}, or \ydiagram[*(red!50)]{1}.
	\item Add to $\sigma$ (or remove from $\tau$) \ydiagram[*(red!50)]{1+1}*[*(white)]{1}, or \ydiagram[*(red!50)]{1,0}*[*(white)]{0,1}.
\end{enumeratei}
It is clear now that $\Da$ is a move-minimizing game digraph, so the (uncolored versions of the) four move-minimizing game questions from Section \ref{s:games} are now on the table.  

\begin{figure}[h]
\resizebox{\textwidth}{!}{
%
%
%
%

\begin{tikzpicture}
\matrix (m) [matrix of math nodes, row sep=1.5em,column sep=1em]
{
 & \TikzPartTwoThree{3,3} &&&&& \ydiagram[*(white)]{0,2+1}*[*(red!50)]{2+1,1+1}*[*(white)\bullet]{1+1,1}*[*(red!50)\bullet]{1,0} &&\\
 & \TikzPartTwoThree{3,2} &&&&& \ydiagram[*(white)]{0,2+1}*[*(red!50)]{0,1+1}*[*(white)\bullet]{1+1,1}*[*(red!50)\bullet]{1,0}*[*(red!50)\bullet]{2+1,0} &&\\
 \TikzPartTwoThree{2,2} && \TikzPartTwoThree{3,1} &&& \ydiagram[*(white)\bullet]{0,2+1}*[*(white)\bullet]{1+1,1}*[*(red!50)\bullet]{1,1+1}*[*(red!50)\bullet]{2+1,0} && \ydiagram[*(white)]{1+1,2+1}*[*(red!50)]{2+1,1+1}*[*(white)\bullet]{0,1}*[*(red!50)\bullet]{1,0} & \\
 & \TikzPartTwoThree{2,1} && \TikzPartTwoThree{3} & \Longleftrightarrow && \ydiagram[*(white)]{0,2+1}*[*(red!50)]{2+1,0}*[*(white)\bullet]{1+1,1}*[*(red!50)\bullet]{1,1+1} && \ydiagram[*(white)]{1+1,2+1}*[*(red!50)]{2+1,1+1}*[*(white)]{0,1}*[*(red!50)]{1,0}\\
 \TikzPartTwoThree{1,1} && \TikzPartTwoThree{2} &&& \ydiagram[*(white)]{0,2+1}*[*(red!50)\bullet]{2+1,0}*[*(white)\bullet]{1+1,1}*[*(red!50)\bullet]{1,1+1} && \ydiagram[*(white)]{0,2+1}*[*(red!50)]{2+1,1+1}*[*(white)\bullet]{1+1,0}*[*(red!50)\bullet]{1,0}*[*(white)]{0,1} & \\
 & \TikzPartTwoThree{1} &&&&& \ydiagram[*(white)]{0,2+1}*[*(red!50)]{0,1+1}*[*(white)\bullet]{1+1,0}*[*(red!50)\bullet]{1,0}*[*(white)]{0,1}*[*(red!50)\bullet]{2+1,0} &&\\
 & \TikzPartTwoThree{0} &&&&& \ydiagram[*(white)]{1+1,2+1}*[*(red!50)]{2+1,1+1}*[*(red!50)\bullet]{1,0}*[*(white)]{0,1} &&\\
};

\path[stealth-]
(m-6-2) edge (m-7-2)
(m-5-1) edge (m-6-2)
(m-5-3) edge (m-6-2)
(m-4-2) edge (m-5-3) edge (m-5-1)
(m-4-4) edge (m-5-3)
(m-3-1) edge (m-4-2)
(m-3-3) edge (m-4-2) edge (m-4-4)
(m-2-2) edge (m-3-3) edge (m-3-1)
(m-1-2) edge (m-2-2)

(m-6-7) edge (m-7-7)
(m-5-6) edge (m-6-7)
(m-5-8) edge (m-6-7)
(m-4-7) edge (m-5-8) edge (m-5-6)
(m-4-9) edge (m-5-8)
(m-3-6) edge (m-4-7)
(m-3-8) edge (m-4-7) edge (m-4-9)
(m-2-7) edge (m-3-8) edge (m-3-6)
(m-1-7) edge (m-2-7);

\end{tikzpicture}

	%
	%
	%
%
}
\caption{$\mathcal{L}(2,3)\cong\mathcal{D}(2,3)$.}
\label{fig:l23d23}
\end{figure}
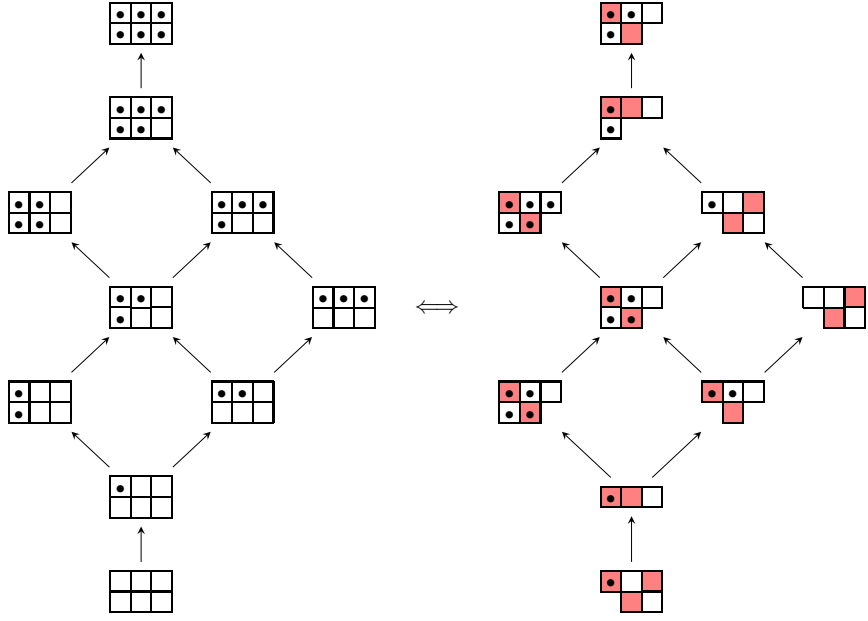

It is clear that both $\Lp$ and $\D$ have the structure of digraphs on the same set of partitions; but similarities in their edge sets are not as clear.  As you can see in Figure \ref{fig:l23d23}, $\D(2,3)\cong \L(2,3)$, but the order of the partitions has changed.  We will show that this is the case in general:  $\La$ and $\Da$ are isomorphic directed graphs, under a permutation on the vertex set that preserves edges.  That is to say, by Definition \ref{d:lA}, that $\D$ has the structure of a distributive lattice where the empty, minimum partition in $\L$ corresponds to a non-empty, minimum partition in $\D$, and where the full, maximum partition in $\L$ corresponds to a non-full, maximum partition in $\Da$.


\paragraph{Domino coordinatizations.}  In order to make our connection between Domino Game digraphs and type $\mathsf{A}$ fundamental lattices precise, we introduce four coordinatizations of $\D$ to parallel the ``part'', ``tab'', ``circ'', and ``diagonal'' coordinatizations of $\L$.  And, as above, we color all the edges with the set $[N-1]$.  The main purpose of the various coordinatizations is that these give the natural progression of the correspondences provided in our main theorem, that $\La\cong\Da$ as edge-colored digraphs. This isomorphism can be shown convincingly through pictures, as in Example \ref{ex:iso}.  Also see Example \ref{ex:big24}, for isomorphic lattices $L_\mathsf{A}(2,4)$ and $\D(2,4)$ including all of the coordinatizations.

\begin{defn}\label{def:dpart} $\Dpart=\left\{k\times(N-k) \text{ partitions} \right \}$.
\begin{itemize}
	\item \textsc{Elements:} $\Dp\subseteq \BZ^k$, and $\sigma=(\sigma_1,\sigma_2,\ldots,\sigma_k)\in\BZ^k$ is in $\Dp$ $\Longleftrightarrow$ $N-k\geq \sigma_1\geq\sigma_2\geq\ldots\geq\sigma_k\geq 0$.
	\item \textsc{Edges:} For $\sigma,\tau\in\Dp$, we have $\sigma\rightarrow\tau$ iff $\tau-\sigma=\beta^\text{part}_l$ for some $l$ with $1\leq l\leq N-1$ where
	
	$\displaystyle{\beta^\text{part}_l=\begin{cases} 
			-2\varepsilon_j & \text{ if } \sigma_j-j = 2l-k, 1 \leq l< N/2 \\
			-\varepsilon_j-\varepsilon_{j+1}& \text{ if } \sigma_{j+1}-j = 2l-k, 1\leq l< N/2 \\
			-\varepsilon_1 & \text{ if } \sigma_1=N-k, l= N/2 \\
			2\varepsilon_j & \text{ if }  \sigma_j-j = 2N-k-2l-1,  N/2 < l \leq N-1 \\
			\varepsilon_j+\varepsilon_{j+1}& \text{ if } \sigma_{j+1}-j = 2N-k-2l,  N/2 < l \leq N-1
	 \end{cases}}$
	
	if $N$ is even, and
	
	$\displaystyle{\beta^\text{part}_l =\begin{cases}
	
			-2\varepsilon_j & \text{ if } \sigma_j-j = 2l-k+1, 1 \leq l < \left\lfloor N/2 \right \rfloor \\
			-\varepsilon_j-\varepsilon_{j+1}& \text{ if } \sigma_{j+1}-j = 2l-k+1, 1 \leq l < \left\lfloor N/2 \right \rfloor \\
			-\varepsilon_1 & \text{ if } \sigma_1=N-k, l=\left\lfloor N/2 \right \rfloor  \\
			2\varepsilon_j & \text{ if }  \sigma_j-j = 2N-k-2l-2,  \left\lceil N/2 \right \rceil \leq l \leq N-1 \\
			\varepsilon_j+\varepsilon_{j+1}& \text{ if } \sigma_{j+1}-j = 2N-k-2l-1,  \left\lceil\! N/2 \!\right \rceil \leq l \leq N-1
		\end{cases}}$
		
		if $N$ is odd.
		
	\item \textsc{Edge colors:} Suppose $\sigma\rightarrow \tau$ in $\Dp$, so $\tau-\sigma=\beta^\text{part}_l$ for some $l$ with $1\leq l\leq N-1$. Then $s\xrightarrow{i}t$ iff $i=l$.

	\item \textsc{Minimum and maximum elements:}  For $N$ even, the minimum element $\mu = (\mu_1, \mu_2, \ldots, \mu_k)$ satisfies $\mu_i = \min\{k-i, N-k\}$ for all $i$, while the maximum element $M = (M_1, M_2, \ldots, M_k)$ satisfies $M_i = \min\{k-i+1, N-k\}$, for all $i$. This minimum can be viewed as shading any box that lies above the diagonal that proceeds from the lower left corner towards the upper right, while for the maximum all boxes on or above the diagonal are shaded.  In the case of $N$ odd, the conditions are reversed: $\mu_i = \min\{k-i+1,N-k\}$ is the minimum element while $M_i = \min\{k-i,N-k\}$ for all $i$ in the maximum element and we can reverse the particular views given.
	
\end{itemize}
\end{defn}

\begin{defn}\label{def:dtab} $\Dtab=\left\{k-\text{element subsets of } [N]\right\}$.
\begin{itemize}
	\item \textsc{Elements:} $\Dt\subseteq\mathcal{P}([N])$, and $S\in\mathcal{P}([N])$ is in $\Dt$ $\Longleftrightarrow$ $|S|=k$.  For $S\in\Dt$, it will be our convention to list the elements of the set in decreasing order.  That is, for $S=\{S_1,S_2,\ldots,S_k\}$ we have $S_i>S_{i+1}$ for $i=1,\ldots, k-1$.
	\item \textsc{Edges:} For $S, T \in \Ltab$, we have $S \rightarrow T$ if and only if $S = (T-\{x\})\cup \{y\}$ where, if $N$ is even,
		
		$\displaystyle{(x,y)=\begin{cases} 
			(2i-1, 2i+1) & \text{ if } 1 \leq i< N/2 \\
			(2i-1, 2i) & \text{ if } i= N/2 \\
			(2N-2i+2, 2N-2i) & \text{ if }  N/2 \leq i \leq N-1
	 \end{cases}}$,
	
		and if $N$ is odd:
		
		$\displaystyle{(x,y) =\begin{cases}
			(2i, 2i+2) & \text{ if } 1 \leq i < \left\lfloor N/2 \right \rfloor \\
			(2i, 2i+1) & \text{ if } i= \left\lfloor N/2 \right \rfloor \\
			(2N-2i+1, 2N-2i-1) & \text{ if }  \left\lceil N/2 \right \rceil \leq i \leq N-1
		\end{cases}}$.
	\item \textsc{Edge colors:} For distinct subsets $S$ and $T$, we say $S\xrightarrow{i}T$ iff $S=(T-\{x\})\cup \{y\}$ for $(x,y)$ as defined above.
	\item \textsc{Visualization:} A $k \times 1$ tableau listing, from greatest to least.  The tableau entries are found by attaching a fully shaded upside down staircase to the left of the associated partition and counting the number of shaded boxes in each row.  See Figure \ref{fig:DTabDiagrams}. An edge of color $i$ exists between two tabular coordinatizations iff the value $x$ becomes $y$ from one diagram to the next, where $x$ and $y$ are defined as above.

		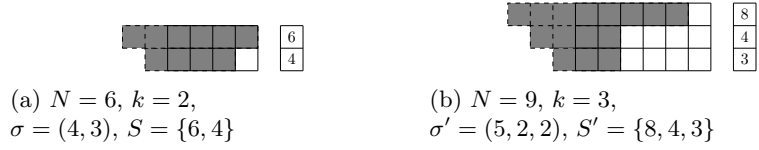
\begin{figure}[H]
		\centering
		\begin{subfigure}[b]{0.45\textwidth}
		\begin{center}
		\scalebox{.6}{
		\begin{tikzpicture}
	\filldraw[fill=gray, dashed] (-.5,0) rectangle (1.5,.5);
	\filldraw[fill=gray, dashed] (-1,.5) rectangle (2,1);

	\foreach \i in {0,...,2}
		\draw (0,\i/2) -- (2,\i/2);
	\foreach \j in {0,...,4}
		\draw (\j/2,0) -- (\j/2,1);
	\draw[dashed] (-.5,.5) -- (-.5,1);
	
	\draw (2.5,0) -- (3,0);
	\draw (2.5,.5) -- (3,.5);
	\draw (2.5,1) -- (3,1);
	\draw (2.5,0) -- (2.5,1);
	\draw (3,0) --(3,1);
	
	\node at (2.75,.25){4};
	\node at (2.75,.75){6};
	
\end{tikzpicture}
		}
		\caption{$N=6$, $k=2$,\\ $\gs = (4,3)$, $S=\{6,4\}$}
		\end{center}
		\end{subfigure}
		\begin{subfigure}[b]{0.45\textwidth}
		\begin{center}
		\scalebox{.6}{
		\begin{tikzpicture}
	\filldraw[fill=gray, dashed] (-1.5,1) rectangle (2.5,1.5);
	\filldraw[fill=gray, dashed] (-1,.5) rectangle (1,1);
	\filldraw[fill=gray,dashed] (-.5,0) rectangle (1,.5);

	\foreach \i in {0,...,3}
		\draw (0,\i/2) -- (3,\i/2);
	\foreach \j in {0,...,6}
		\draw (\j/2,0) -- (\j/2,1.5);
	\draw[dashed] (-.5,.5) -- (-.5,1.5);
	\draw[dashed] (-1,1) -- (-1,1.5);

	\draw (3.5,0) -- (4,0);
	\draw (3.5,.5) -- (4,.5);
	\draw (3.5,1) -- (4,1);
	\draw (3.5,1.5) -- (4,1.5);
	\draw (3.5,0) -- (3.5,1.5);
	\draw (4,0) --(4,1.5);
	
	\node at (3.75,.25){3};
	\node at (3.75,.75){4};
	\node at (3.75,1.25){8};
	
\end{tikzpicture}
		}
		\caption{$N=9$, $k=3$, \\ $\gs' = (5,2,2)$, $S'=\{8,4,3\}$}
		\end{center}
		\end{subfigure}
		\caption{\label{fig:DTabDiagrams}%
		Two partitions with augmented staircases and their tableaux.}
		\end{figure}
\end{itemize}
\end{defn}

\begin{defn}\label{def:dcirc} $\Dcirc=\left\{\text{Length $N$ binary sequences with $k$ 1's}\right\}$.
\begin{itemize}
	\item \textsc{Elements:} $\Dc\subseteq\BZ_2^N$, and $s=(s_1,s_2,\ldots,s_N)\in\BZ_2^{N}$ $\Longleftrightarrow$ $\sum s_i=k$.
	\item \textsc{Edges:} For $s,t\in\Dc$, we have $s\rightarrow t$ iff $t-s= \beta^\text{circ}_l$ for some $l$ with $1\leq l\leq N-1$ where
	
		$\displaystyle{\beta^\text{circ}_l=\begin{cases} 
			\varepsilon_{2l-1}+ \varepsilon_{2l+1}& \text{ if } 1 \leq l< N/2 \\
			\varepsilon_{2l-1}-\varepsilon_{2l}& \text{ if } l= N/2 \\
			\varepsilon_{2N-2l+2}-\varepsilon_{2N-2l}& \text{ if }  N/2 \leq l \leq N-1
	 \end{cases}}$, if $N$ is even.
	
	$\displaystyle{\beta^\text{circ}_l =\begin{cases}
			\varepsilon_{2l}+ \varepsilon_{2l+2}& \text{ if } 1 \leq l \leq \left\lfloor N/2 \right \rfloor \\
			\varepsilon_{2l}-\varepsilon_{2l+1}& \text{ if } l= \left\lfloor N/2 \right \rfloor \\
			\varepsilon_{2N-2l+1}- \varepsilon_{2N-2l-1}& \text{ if }  \left\lceil N/2 \right \rceil \leq l \leq N-1
		\end{cases}}$, if $N$ is odd.
	\item \textsc{Edge colors:} Suppose $s\rightarrow t$ in $\Dc$, so $t-s=\beta^\text{diag}_l$ for some $l$ with $1\leq l\leq N-1$. Then $s\xrightarrow{i}t$ iff $i=l$.
	\item \textsc{Visualization:}  A $2 \times \left\lceil N/2 \right\rceil$ tableau (with the two rightmost boxes joined into one tall box when $N$ is odd) numbered as in \ref{fig:DCircleDiagrams}: we number the boxes top to bottom, from left to right, beginning in the top left box. The box labeled $i$ is marked with a dot iff $s_i=1$, $1 \leq i\leq N$.  Note the identification of edge colors is independent of the box numbers; it depends only on the underlying shape of the circle diagram.  That is, the identification of edge colors is identical in both $\Lc$ and $\Dc$.  
	
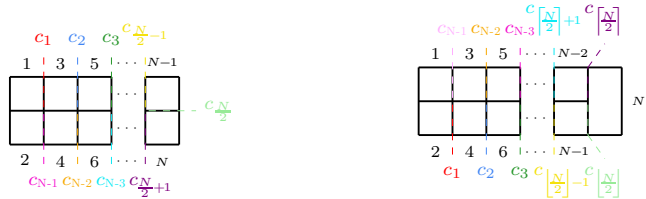
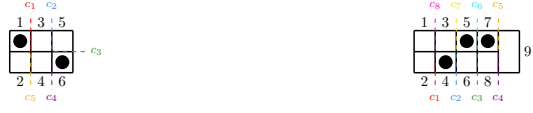
\begin{figure}[H]
\centering
\begin{subfigure}[b]{0.45\textwidth}
\begin{center}
\scalebox{.9}{%
\begin{tikzpicture}

\foreach \place/\name in {{(0,.5)/ta}, {(.5,.5)/tb}, {(1,.5)/tc}, {(1.5,.5)/td}, {(2,.5)/te}, {(2.5,.5)/tf}, {(0,-.5)/ba}, {(.5,-.5)/bb}, {(1,-.5)/bc}, {(1.5,-.5)/bd}, {(2,-.5)/be}, {(2.5,-.5)/bf}, {(0,0)/ma},  {(1.5,0)/md}, {(2,0)/me}, {(2.5,0)/mf}}
    \node[circle, fill=white] (\name) at \place {};

\draw[line width=.8pt] (ta.base) -- (td.base);
\draw[line width=.8pt] (te.base) -- (tf.base);
\draw[line width=.8pt] (ma.base) -- (md.base);
\draw[line width=.8pt] (me.base) -- (mf.base);
\draw[line width=.8pt] (ba.base) -- (bd.base);
\draw[line width=.8pt] (be.base) -- (bf.base);
\draw[line width=.8pt] (ta.base) -- (ba.base);
\draw[line width=.8pt] (tb.base) -- (bb.base);
\draw[line width=.8pt] (tc.base) -- (bc.base);
\draw[line width=.8pt] (td.base) -- (bd.base);
\draw[line width=.8pt] (te.base) -- (be.base);
\draw[line width=.8pt] (tf.base) -- (bf.base);

\node [font=\scriptsize] at (0.25,.7){1};
\node [font=\scriptsize] at (.75,.7){3};
\node [font=\scriptsize] at (1.25,.7){5};
\node [font=\tiny] at (1.75,.7){$\cdots$};
\node [font=\tiny] at (2.25,.7){$N\!\!-\!\!1$};

\node [font=\scriptsize] at (0.22,-.73){2};
\node [font=\scriptsize] at (.75,-.75){4};
\node [font=\scriptsize] at (1.25,-.75){6};
\node [font=\tiny] at (1.75,-.75){$\cdots$};
\node [font=\tiny] at (2.25,-.75){$N$};

\node [font=\tiny] at (1.75,-0.25){$\cdots$};
\node [font=\tiny] at (1.75,0.25){$\cdots$};

\draw[red, dashed] (.5,.05) -- (.5,.85) node[above, font=\footnotesize]{$c_1$};
\draw[blue, dashed] (1,.05) -- (1,.85) node[above, font=\footnotesize]{$c_2$};
\draw[green, dashed] (1.5,.05) -- (1.5,.85) node[above, font=\footnotesize]{$c_3$};
\draw[gold,dashed] (2,.05) -- (2,.85) node[above, font=\footnotesize]{$c_\subsc{$\frac{N}{2}\!\!-\!\!1$}$};
\draw[seafoam, dashed] (2.05,0) -- (2.75,0) node[right, font=\footnotesize]{$c_\subsc{$\frac{N}{2}$}$};
\draw[violet, dashed] (2,-.05) -- (2,-.85) node[below, font=\footnotesize]{$\;\;c_\subsc{$\!\frac{N}{2}\!\!+\!\!1$}$};
\draw[cyan, dashed] (1.5,-.05) -- (1.5,-.85) node[below, font=\footnotesize]{$c_\subsc{N-3}$};
\draw[orange, dashed] (1,-.05) -- (1,-.85) node[below, font=\footnotesize]{$c_\subsc{N-2}$};
\draw[magenta, dashed] (.5,-.05) -- (.5,-.85) node[below, font=\footnotesize]{$c_\subsc{N-1}$};

\end{tikzpicture}%
}
\end{center}
\caption{Circle diagram for $N$ even}
\end{subfigure}
\begin{subfigure}[b]{0.45\textwidth}
\begin{center}
\scalebox{.9}{%
\begin{tikzpicture}

\foreach \i in {0,...,6}
	\draw[line width=.8pt] (\i/2,0) -- (\i/2,1);

\foreach \j in {0,...,2}
	\draw[line width=.8pt] (0,\j/2) -- (1.5,\j/2);
\draw[line width=.8pt] (2,0) -- (3,0);
\draw[line width=.8pt] (2,.5) -- (2.5,.5);
\draw[line width=.8pt] (2,1) -- (3,1);

\draw[red, dashed] (.5,.45) -- (.5,-.35) node[below, font=\footnotesize]{$c_1$};
\draw[blue, dashed] (1,.45) -- (1,-.35) node[below, font=\footnotesize]{$c_2$};
\draw[green, dashed] (1.5,.45) -- (1.5,-.35) node[below, font=\footnotesize]{$c_3$};
\draw[gold,dashed] (2,.45) -- (2,-.35) node[below, font=\footnotesize]{$\;\;\;c_\subsc{$\left\lfloor\!\!\frac{N}{2}\!\!\right\rfloor\!\!-\!\!1$}$};
\draw[seafoam, dashed] (2.5,.45) -- (2.5,0);
\draw[seafoam, dashed] (2.5,0) -- (2.75,-.35) node[below, font=\footnotesize]{$c_\subsc{$\left\lfloor\!\!\frac{N}{2}\!\!\right\rfloor$}$};
\draw[cyan, dashed] (2,.55) -- (2,1.35) node[above, font=\footnotesize]{$c_\subsc{$\left\lceil\!\!\frac{N}{2}\!\!\right\rceil\!\!+\!\!1$}$};
\draw[violet, dashed] (2.5,.55) -- (2.5,1);
\draw[violet, dashed] (2.5,1) -- (2.75,1.35) node[above, font=\footnotesize]{$c_\subsc{$\left\lceil\!\!\frac{N}{2}\!\!\right\rceil$}$};
\draw[magenta, dashed] (1.5,.55) -- (1.5,1.35) node[above, font=\footnotesize]{$c_\subsc{N-3}$};
\draw[orange, dashed] (1,.55) -- (1,1.35) node[above, font=\footnotesize]{$c_\subsc{N-2}$};
\draw[lavender, dashed] (.5,.55) -- (.5,1.35) node[above, font=\footnotesize]{$c_\subsc{N-1}$};

\node[font=\scriptsize] at (0.25,1.2){$1$};
\node[font=\scriptsize] at (.75,1.2){$3$};
\node[font=\scriptsize] at (1.25,1.2){$5$};
\node[font=\tiny] at (1.75,1.2){$\cdots$};
\node[font=\tiny] at (2.2,1.2){$\;\;N\!\!-\!\!2$};
\node[font=\tiny] at (3.25,.5){$N$};

\node[font=\scriptsize] at (0.25,-.25){$2$};
\node[font=\scriptsize] at (.75,-.25){$4$};
\node[font=\scriptsize] at (1.25,-.25){$6$};
\node[font=\tiny] at (1.75,-.25){$\cdots$};
\node[font=\tiny] at (2.25,-.25){$\;N\!\!-\!\!1$};

\node[font=\tiny] at (1.75,0.25){$\cdots$};
\node[font=\tiny] at (1.75,0.75){$\cdots$};

\end{tikzpicture}
}
\end{center}
\caption{Circle diagram for $N$ odd}
\end{subfigure}

\caption{\label{fig:DCircleDiagrams}
Generic circle diagram labellings.}
\end{figure}

\begin{figure}[H]
\centering
\begin{subfigure}[b]{0.45\textwidth}
\begin{center}
\scalebox{.56}{%
\begin{tikzpicture}

\foreach \place/\name in {{(0,.5)/ta}, {(.5,.5)/tb}, {(1,.5)/tc}, {(1.5,.5)/td}, {(0,-.5)/ba}, {(.5,-.5)/bb}, {(1,-.5)/bc}, {(1.5,-.5)/bd}, {(0,0)/ma},  {(1.5,0)/md}, {(1.5,0)/md}}
    \node[circle, fill=white] (\name) at \place {};

\node[circle, fill=black] (dot1) at (.25,.25){};
\node[circle, fill=black] (dot2) at (1.25,-.25){};

\draw[line width=1pt] (ta.base) -- (td.base);
\draw[line width=1pt]  (ma.base) -- (md.base);
\draw[line width=1pt]  (ba.base) -- (bd.base);
\draw[line width=1pt]  (ta.base) -- (ba.base);
\draw[line width=1pt]  (tb.base) -- (bb.base);
\draw[line width=1pt]  (tc.base) -- (bc.base);
\draw[line width=1pt]  (td.base) -- (bd.base);

\draw[dashed,red] (.5,.05) -- (.5,.9)node[above, font=\footnotesize]{$c_1$};
\draw[dashed,blue] (1,.05) -- (1,.9)node[above, font=\footnotesize]{$c_2$};
\draw[dashed,green] (1.05,0) -- (1.8,0)node[right, font=\footnotesize]{$c_3$};
\draw[dashed,violet] (1,-.05) -- (1,-.9)node[below, font=\footnotesize]{$c_4$};
\draw[dashed,orange] (.5,-.05) -- (.5,-.9)node[below, font=\footnotesize]{$c_5$};

\node (v1) at (0.25,.7){1};
\node (v2) at (.75,.7){3};
\node (v3) at (1.25,.7){5};

\node (v6) at (0.25,-.7){2};
\node (v7) at (.75,-.7){4};
\node (v8) at (1.25,-.7){6};

\end{tikzpicture}%
}
\end{center}
\caption{$N=6$, $k=2$,\\ $\gs = (4,3)$, $S=\{1,6\}$}
\end{subfigure}
\begin{subfigure}[b]{0.45\textwidth}
\begin{center}
\scalebox{.56}{%

\begin{tikzpicture}

\foreach \place/\name in {{(0,.5)/ta}, {(.5,.5)/tb}, {(1,.5)/tc}, {(1.5,.5)/td}, {(2,.5)/te}, {(2.5,.5)/tf}, {(0,-.5)/ba}, {(.5,-.5)/bb}, {(1,-.5)/bc}, {(1.5,-.5)/bd}, {(2,-.5)/be}, {(2.5,-.5)/bf}, {(0,0)/ma}, {(2,0)/me}}
    \node[circle, fill=white] (\name) at \place {};

\node[fill=white] (tl) at (2.1,.5) {};
\node[fill=white] (tr) at (2.6,.5) {};
\node[fill=white] (bl) at (2.1,0) {};
\node[fill=white] (br) at (2.6,0) {};

\node[circle, fill=black] (dot1) at (.75,-.25) {};
\node[circle, fill=black] (dot2) at (1.25,.25){};
\node[circle, fill=black] (dot3) at (1.75,.25){};

\draw[line width=1pt] (ta.base) -- (tf.base);
\draw[line width=1pt] (ma.base) -- (me.base);
\draw[line width=1pt] (ba.base) -- (bf.base);
\draw[line width=1pt] (ta.base) -- (ba.base);
\draw[line width=1pt] (tb.base) -- (bb.base);
\draw[line width=1pt] (tc.base) -- (bc.base);
\draw[line width=1pt] (td.base) -- (bd.base);
\draw[line width=1pt] (te.base) -- (be.base);
\draw[line width=1pt] (tf.base) -- (bf.base);

\draw[dashed,red] (.5,-.05) -- (.5,-.9)node[below, font=\footnotesize]{$c_1$};
\draw[dashed,blue] (1,-.05) -- (1,-.9)node[below, font=\footnotesize]{$c_2$};
\draw[dashed,green] (1.5,-.05) -- (1.5,-.9)node[below, font=\footnotesize]{$c_3$};
\draw[dashed,violet] (2,-.05) -- (2,-.9)node[below, font=\footnotesize]{$c_4$};
\draw[dashed,orange] (2,.05) -- (2,.9)node[above, font=\footnotesize]{$c_5$};
\draw[dashed,cyan] (1.5,.05) -- (1.5,.9)node[above, font=\footnotesize]{$c_6$};
\draw[dashed,gold] (1,.05) -- (1,.9)node[above, font=\footnotesize]{$c_7$};
\draw[dashed,magenta] (.5,.05) -- (.5,.9)node[above, font=\footnotesize]{$c_8$};

\node (v1) at (0.25,.7){1};
\node (v2) at (.75,.7){3};
\node (v3) at (1.25,.7){5};
\node (v4) at (1.75,.7){7};
\node (v5) at (2.7,0){9};

\node (v6) at (0.25,-.7){2};
\node (v7) at (.75,-.7){4};
\node (v8) at (1.25,-.7){6};
\node (v9) at (1.75,-.7){8};

\end{tikzpicture}

}
\end{center}
\caption{$N=9$, $k=3$, \\ $\gs' = (5,2,2)$, $S'=\{4,5,7\}$}
\end{subfigure}

\caption{
Two circle diagrams for the given tabular coordinatizations.}

\end{figure}
\end{itemize}
\end{defn}

\begin{defn}\label{def:ddiag} $\Ddiag=\left\{\text{diagonal coords. of $k\!\times\!(N\!-\!k)$ partitions}\right\}$.
\begin{itemize}
	\item \textsc{Elements:} $\Dd\subseteq\BZ^{N-1}$, and $\ds=(\ds_1,\ds_2,\ldots,\ds_{N-1})\in\Dd$ iff
	\begin{itemize}
		\item $\ds_i\leq \min\{i,k\}$, for $1\leq i\leq N-k$,
		\item $\ds_i\leq \min\{N-i, N-k\}$, for $N-k\leq i\leq N-1$, 
		\item For $1\leq i\leq N-k-1$, $\ds_{i+1}=\ds_i$ or $\ds_i+1$, and
		\item For $N-k\leq i\leq N-2$, $\ds_{i+1}=\ds_i$ or $\ds_i-1$.
	\end{itemize}
	\item \textsc{Edges:} For $\ds,\dt\in\Dd$, we have $\ds\rightarrow \dt$, iff $\dt-\ds=\beta^\text{diag}_l$ for some $l$ with $1\leq l\leq N-1$ where
	
		$\displaystyle{\beta^\text{diag}_l=\begin{cases} 
				-\varepsilon_{N-2l}- \varepsilon_{N-2l+1}& \text{ if } 1\leq l< N/2 \\
				-\varepsilon_1& \text{ if } l= N/2 \\
				\varepsilon_{2l-N}+ \varepsilon_{2l-N-1}& \text{ if }  N/2 <l \leq N-1
		 \end{cases}}$, if $N$ is even.
	
		$\displaystyle{\beta^\text{diag}_l =\begin{cases}
				-\varepsilon_{N-2l}- \varepsilon_{N-2l-1}& \text{ if } 1 \leq l \leq \left\lfloor N/2 \right \rfloor \\
				-\varepsilon_{1}& \text{ if } l= \left\lfloor N/2 \right \rfloor \\
				\varepsilon_{2l-N}+ \varepsilon_{2l-N+1}& \text{ if }  \left\lceil N/2 \right \rceil \leq l \leq N-1
			\end{cases}}$, if $N$ is odd.
	\item \textsc{Edge colors}: Suppose $\ds\rightarrow \dt$ in $\Dd$, so $\dt-\ds=\beta^\text{diag}_l$ for some $l$ with $1\leq l\leq N-1$. Then $\ds\xrightarrow{i}\dt$ iff $i=l$.
\end{itemize}
\end{defn}

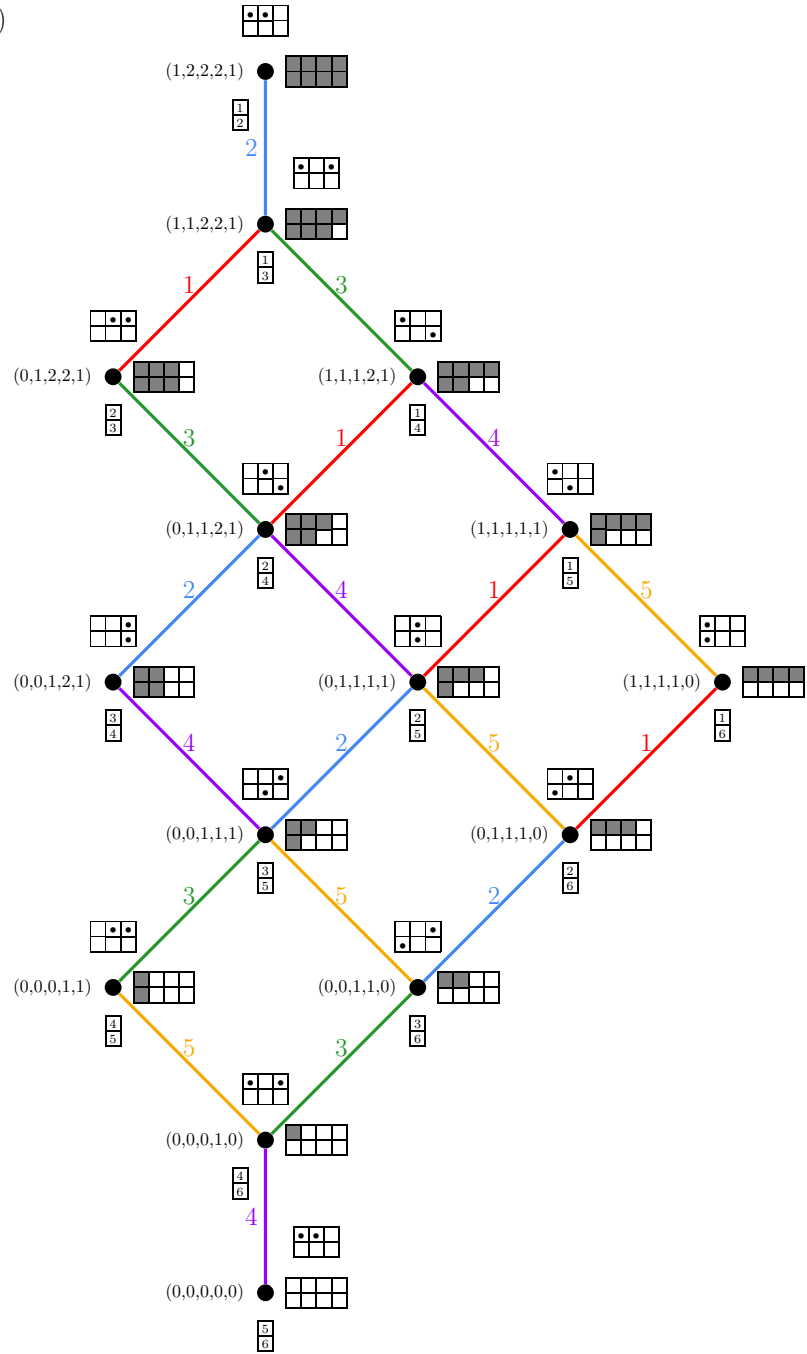
\begin{figure}[!p]
\begin{center}
\resizebox{\textwidth}{.935\textheight}{%
\begin{tikzpicture}

\foreach \place/\name in {{(0,0)/L1c}, {(0,3)/L2c}, {(-3,6)/L3l}, {(3,6)/L3r}, {(0,9)/L4c}, {(6,9)/L4rr}, {(-3,12)/L5l}, {(3,12)/L5r}, {(9,12)/L5rrr}, {(0,15)/L6c}, {(6,15)/L6rr}, {(-3,18)/L7l}, {(3,18)/L7r}, {(0,21)/L8c}, {(0,24)/L9c}}
    \node[circle, fill=black] (\name) at \place {};
\foreach \source/\dest in {L4rr/L5rrr, L5r/L6rr, L6c/L7r, L7l/L8c}
		\draw[line width=1.8pt, red] (\source) -- node[above, font=\Large] {1} (\dest);
\foreach \source/\dest in {L3r/L4rr, L4c/L5r, L5l/L6c}
    \draw[line width=1.8pt, blue] (\source) -- node[above, font=\Large] {2} (\dest);
		\draw[line width=1.8pt, blue] (L8c) -- node[left, font=\Large]{2} (L9c);
\foreach \source/\dest in {L2c/L3r, L3l/L4c, L6c/L7l, L7r/L8c}
    \draw[line width=1.8pt, green] (\source) -- node[above, font=\Large] {3} (\dest);
\foreach \source/\dest in {L4c/L5l, L5r/L6c, L6rr/L7r}
		\draw[line width=1.8pt, purple] (\source) -- node[above, font=\Large] {4} (\dest);
		\draw[line width=1.8pt, purple] (L1c) -- node[left, font=\Large] {4} (L2c);
\foreach \source/\dest in {L2c/L3l, L3r/L4c, L4rr/L5r, L5rrr/L6rr}
    \draw[line width=1.8pt, orange] (\source) -- node[above, font=\Large] {5} (\dest);

\node (L1cCirc) at (1,1) {\circdiag{0,2}};
\node (L1cTab) at (0,-.85) {\tab{5}{6}};
\node (L1cPart) at (1,0) {\partpic{0,0}};
\node (L1cDiag) at (-1.2,0) {(0,0,0,0,0)};

\node (L2cCirc) at (0,4) {\ydiagram[*(white)\bullet]{0,1}*[*(white)\bullet]{0,2+1}*[*(white)]{3,3}};
\node (L2cTab) at (-.5,2.15) {\tab{4}{6}};
\node (L2cPart) at (1,3) {\partpic{1,0}};
\node (L2cDiag) at (-1.2,3) {(0,0,0,1,0)};

\node (L3lCirc) at (-3,7) {\circdiag{0,1+2}};
\node (L3lTab) at (-3,5.15) {\tab{4}{5}};
\node (L3lPart) at (-2,6) {\partpic{1,1}};
\node (L3lDiag) at (-4.2,6) {(0,0,0,1,1)};

\node (L3rCirc) at (3,7) {\circdiag{2+1,1}};
\node (L3rTab) at (3,5.15) {\tab{3}{6}};
\node (L3rPart) at (4,6) {\partpic{2,0}};
\node (L3rDiag) at (1.8,6) {(0,0,1,1,0)};

\node (L4cCirc) at (0,10) {\circdiag{2+1,1+1}};
\node (L4cTab) at (0,8.15) {\tab{3}{5}};
\node (L4cPart) at (1,9) {\partpic{2,1}};
\node (L4cDiag) at (-1.2,9) {(0,0,1,1,1)};

\node (L4rCirc) at (6,10) {\circdiag{1+1,1}};
\node (L4rTab) at (6,8.15) {\tab{2}{6}};
\node (L4rPart) at (7,9) {\partpic{3,0}};
\node (L4rDiag) at (4.8,9) {(0,1,1,1,0)};

\node (L5lCirc) at (-3,13) {\circdiag{2+1,2+1}};
\node (L5lTab) at (-3,11.15) {\tab{3}{4}};
\node (L5lPart) at (-2,12) {\partpic{2,2}};
\node (L5lDiag) at (-4.2,12) {(0,0,1,2,1)};

\node (L5rCirc) at (3,13) {\circdiag{1+1,1+1}};
\node (L5rTab) at (3,11.15) {\tab{2}{5}};
\node (L5rPart) at (4,12) {\partpic{3,1}};
\node (L5rDiag) at (1.8,12) {(0,1,1,1,1)};

\node (L5rrrCirc) at (9,13) {\circdiag{1,1}};
\node (L5rrrTab) at (9,11.15) {\tab{1}{6}};
\node (L5rrrPart) at (10,12) {\partpic{4,0}};
\node (L5rrrDiag) at (7.8,12) {(1,1,1,1,0)};

\node (L6cCirc) at (0,16) {\circdiag{1+1,2+1}};
\node (L6cTab) at (0,14.15) {\tab{2}{4}};
\node (L6cPart) at (1,15) {\partpic{3,2}};
\node (L6cDiag) at (-1.2,15) {(0,1,1,2,1)};

\node (L6rrCirc) at (6,16) {\circdiag{1,1+1}};
\node (L6rrTab) at (6,14.15) {\tab{1}{5}};
\node (L6rrPart) at (7,15) {\partpic{4,1}};
\node (L6rrDiag) at (4.8,15) {(1,1,1,1,1)};

\node (L7lCirc) at (-3,19) {\circdiag{1+2,0}};
\node (L7lTab) at (-3,17.15) {\tab{2}{3}};
\node (L7lPart) at (-2,18) {\partpic{3,3}};
\node (L7lDiag) at (-4.2,18) {(0,1,2,2,1)};

\node (L7rCirc) at (3,19) {\circdiag{1,2+1}};
\node (L7rTab) at (3,17.15) {\tab{1}{4}};
\node (L7rPart) at (4,18) {\partpic{4,2}};
\node (L7rDiag) at (1.8,18) {(1,1,1,2,1)};

\node (L8cCirc) at (1,22) {\ydiagram[*(white)\bullet]{1,0}*[*(white)\bullet]{2+1,0}*[*(white)]{3,3}};
\node (L8cTab) at (0,20.15) {\tab{1}{3}};
\node (L8cPart) at (1,21) {\partpic{4,3}};
\node (L8cDiag) at (-1.2,21) {(1,1,2,2,1)};

\node (L9cCirc) at (0,25) {\circdiag{2,0}};
\node (L9cTab) at (-.5,23.15) {\tab{1}{2}};
\node (L9cPart) at (1,24) {\partpic{4,4}};
\node (L9cDiag) at (-1.2,24) {(1,2,2,2,1)};

\node[font=\fontsize{16}{18}\selectfont] (name) at (-6,25) {$\L(2,4)$};

\end{tikzpicture} }
\end{center}
\caption{\label{fig:L24Full} $\L(2,4)$ with all four coordinatizations and labeled edges.}
\end{figure}

\begin{Ex}\label{ex:big24}
This example completely describes the coordinatizations in Sections \ref{s:classicalL} and \ref{s:dominos}.  In Figure \ref{fig:L24Full}, we view $\L(2,4)$ and then the corresponding $\D(2,4)$ in Figure \ref{fig:D24Full}. In both cases, each vertex is surrounded by the corresponding coordinatizations: circle diagram, partition, tableau, and diagonal coordinates.  Also each edge is marked with its color value.
\end{Ex}

\begin{figure}[!p]
\begin{center}
\resizebox{\textwidth}{.935\textheight}{%
 \begin{tikzpicture}

\foreach \place/\name in {{(0,0)/L1c}, {(0,3)/L2c}, {(-3,6)/L3l}, {(3,6)/L3r}, {(0,9)/L4c}, {(6,9)/L4rr}, {(-3,12)/L5l}, {(3,12)/L5r}, {(9,12)/L5rrr}, {(0,15)/L6c}, {(6,15)/L6rr}, {(-3,18)/L7l}, {(3,18)/L7r}, {(0,21)/L8c}, {(0,24)/L9c}}
    \node[circle, fill=black] (\name) at \place {};
\foreach \source/\dest in {L4rr/L5rrr, L5r/L6rr, L6c/L7r, L7l/L8c}
		\draw[line width=1.8pt, red] (\source) -- node[above, font=\Large] {1} (\dest);
\foreach \source/\dest in {L3r/L4rr, L4c/L5r, L5l/L6c}
    \draw[line width=1.8pt, blue] (\source) -- node[above, font=\Large] {2} (\dest);
		\draw[line width=1.8pt, blue] (L8c) -- node[left, font=\Large]{2} (L9c);
\foreach \source/\dest in {L2c/L3r, L3l/L4c, L6c/L7l, L7r/L8c}
    \draw[line width=1.8pt, green] (\source) -- node[above, font=\Large] {3} (\dest);
\foreach \source/\dest in {L4c/L5l, L5r/L6c, L6rr/L7r}
		\draw[line width=1.8pt, purple] (\source) -- node[above, font=\Large] {4} (\dest);
		\draw[line width=1.8pt, purple] (L1c) -- node[left, font=\Large] {4} (L2c);
\foreach \source/\dest in {L2c/L3l, L3r/L4c, L4rr/L5r, L5rrr/L6rr}
    \draw[line width=1.8pt, orange] (\source) -- node[above, font=\Large] {5} (\dest);

\node (L1cCirc) at (1,1) {\circdiag{0,2}};
\node (L1cTab) at (0,-.85) {\tab{4}{2}};
\node (L1cPart) at (1,0) {\partpic{2,1}};
\node (L1cDiag) at (-1.2,0) {(0,0,1,1,1)};

\node (L2cCirc) at (0,4) {\ydiagram[*(white)\bullet]{0,1}*[*(white)\bullet]{0,2+1}*[*(white)]{3,3}};
\node (L2cTab) at (-.5,2.15) {\tab{6}{2}};
\node (L2cPart) at (1,3) {\partpic{4,1}};
\node (L2cDiag) at (-1.2,3) {(1,1,1,1,1)};

\node (L3lCirc) at (-3,7) {\circdiag{0,1+2}};
\node (L3lTab) at (-3,5.15) {\tab{6}{4}};
\node (L3lPart) at (-2,6) {\partpic{4,3}};
\node (L3lDiag) at (-4.2,6) {(1,1,2,2,1)};

\node (L3rCirc) at (3,7) {\circdiag{2+1,1}};
\node (L3rTab) at (3,5.15) {\tab{5}{2}};
\node (L3rPart) at (4,6) {\partpic{3,1}};
\node (L3rDiag) at (1.8,6) {(0,1,1,1,1)};

\node (L4cCirc) at (0,10) {\circdiag{2+1,1+1}};
\node (L4cTab) at (0,8.15) {\tab{5}{4}};
\node (L4cPart) at (1,9) {\partpic{3,3}};
\node (L4cDiag) at (-1.2,9) {(0,1,2,2,1)};

\node (L4rCirc) at (6,10) {\circdiag{1+1,1}};
\node (L4rTab) at (6,8.15) {\tab{3}{2}};
\node (L4rPart) at (7,9) {\partpic{1,1}};
\node (L4rDiag) at (4.8,9) {(0,0,0,1,1)};

\node (L5lCirc) at (-3,13) {\circdiag{2+1,2+1}};
\node (L5lTab) at (-3,11.15) {\tab{6}{5}};
\node (L5lPart) at (-2,12) {\partpic{4,4}};
\node (L5lDiag) at (-4.2,12) {(1,2,2,2,1)};

\node (L5rCirc) at (3,13) {\circdiag{1+1,1+1}};
\node (L5rTab) at (3,11.15) {\tab{4}{3}};
\node (L5rPart) at (4,12) {\partpic{2,2}};
\node (L5rDiag) at (1.8,12) {(0,0,1,2,1)};

\node (L5rrrCirc) at (9,13) {\circdiag{1,1}};
\node (L5rrrTab) at (9,11.15) {\tab{2}{1}};
\node (L5rrrPart) at (10,12) {\partpic{0,0}};
\node (L5rrrDiag) at (7.8,12) {(0,0,0,0,0)};

\node (L6cCirc) at (0,16) {\circdiag{1+1,2+1}};
\node (L6cTab) at (0,14.15) {\tab{6}{3}};
\node (L6cPart) at (1,15) {\partpic{4,2}};
\node (L6cDiag) at (-1.2,15) {(1,1,1,2,1)};

\node (L6rrCirc) at (6,16) {\circdiag{1,1+1}};
\node (L6rrTab) at (6,14.15) {\tab{4}{1}};
\node (L6rrPart) at (7,15) {\partpic{2,0}};
\node (L6rrDiag) at (4.8,15) {(0,0,1,1,0)};

\node (L7lCirc) at (-3,19) {\circdiag{1+2,0}};
\node (L7lTab) at (-3,17.15) {\tab{5}{3}};
\node (L7lPart) at (-2,18) {\partpic{3,2}};
\node (L7lDiag) at (-4.2,18) {(0,1,1,2,1)};

\node (L7rCirc) at (3,19) {\circdiag{1,2+1}};
\node (L7rTab) at (3,17.15) {\tab{6}{1}};
\node (L7rPart) at (4,18) {\partpic{4,0}};
\node (L7rDiag) at (1.8,18) {(1,1,1,1,0)};

\node (L8cCirc) at (1,22) {\ydiagram[*(white)\bullet]{1,0}*[*(white)\bullet]{2+1,0}*[*(white)]{3,3}};
\node (L8cTab) at (0,20.15) {\tab{5}{1}};
\node (L8cPart) at (1,21) {\partpic{3,0}};
\node (L8cDiag) at (-1.2,21) {(0,1,1,1,0)};

\node (L9cCirc) at (0,25) {\circdiag{2,0}};
\node (L9cTab) at (-.5,23.15) {\tab{3}{1}};
\node (L9cPart) at (1,24) {\partpic{1,0}};
\node (L9cDiag) at (-1.2,24) {(0,0,0,1,0)};

\node[font=\fontsize{16}{18}\selectfont] (name) at (-6,25) {$\D(2,4)$};

\end{tikzpicture} }
\end{center}
\caption{\label{fig:D24Full} $\D(2,4)$ with all four coordinatizations and labeled edges.}
\end{figure}

\paragraph{The Coordinate Isomorphisms.}  The purpose of this section is to make precise the correspondence $\La\cong\Da$ by defining isomorphisms 
\[\Dd\leftrightarrow\Dp\leftrightarrow\Dt\leftrightarrow\Dc.\]
Ultimately, we will show that $\Lc\cong\Dc$, and therefore all of the coordinatizations of both $\La$ and $\Da$ are isomorphic.  Moreover, each map preserves edges and edge colors, so that it will be accurate to say that $\L\cong\D$ as diamond-colored distributive lattices, and the answers to our move-minimizing questions for $\D$ will follow from this isomorphism.  

These seemingly complicated maps are really not hard to understand pictorially.  Note that we do not include the proofs that each map is indeed an edge color preserving isomorphism here, as the proofs intuitive and the description given in Example \ref{ex:iso} is convincing.

\begin{defn}\label{d:part-tab} Define $\gampt:\Dp\rightarrow\Dt$ as $\gampt(\gs) = \{\gs_j+k-j+1\}_{j=1}^k$ and its inverse $\gamtp:\Dt\rightarrow\Dp$ by $\gamtp(S) = (S_j-k+j-1)_{j=1}^k$ where again, $S=\{S_1>S_2> \cdots > S_k\}$.
\end{defn}

\begin{defn}\label{d:tab-circ} Define $\gamtc:\Dt\rightarrow\Dc$ be defined as $\gamtc(S) = \{s = (s_1, s_2, \ldots, s_N) | s_j = 1 \text{ if } j\in S, 0 \text{ otherwise}\}$ and its inverse $\gamct:\Dc\rightarrow\Dt$ by $\gamct:\Dtab \rightarrow\Dpart$ be defined as $\gamct(s) = \{j|s_j = 1\}_{j=1}^k$.
\end{defn}

\begin{defn}\label{d:part-diag} We define a map $\gampd:\Dp\rightarrow\Dd$ by 
\begin{equation}\label{e:gampd}
\ds_i=\begin{cases}
	\max\{0,j\mid \sigma_j\geq (N-k)-(i-j)\}_{j=1}^{m_i} & \text{ if } 1\leq i\leq N-k\\
	\max\{0,j\mid \sigma_{(j+i)-(N-k)}\geq j\}_{j=1}^{n_i} & \text{ if } N-k \leq i\leq N-1.
	\end{cases}
\end{equation}
The map $\gamdp:\Dd\rightarrow\Dp$ is defined by 
\begin{equation}\label{e:gamdp}
	\sigma_i=\sum_{\substack{j=i \\ \ds_j\geq i}}^{N-k}1 + \sum_{\substack{l=1 \\ \ds_{N-k+l}\geq (i-l)}}^{i-1}1.
\end{equation}
The visualization in Definition \ref{def:ldiag} describes the isomorphism (cf. Definition \ref{def:ldiag}).
\end{defn}

\paragraph{The Main Theorem.}  Note that $\Lc$ and $\Dc$ are exactly equal on the level of circle diagrams, and have obviously equivalent directed edges and colors. It is simply that the numbering of the boxes is different. (See Definitions \ref{def:lcirc} and \ref{def:dcirc}.) So, to define an isomorphism $\Lc\rightarrow\Dc$ we simply define the appropriate permutation $\pi \in S^N$ that renumbers the boxes.  The definition depends on the parity of $N$.

For $N$ even, we have
\begin{equation}\label{e:circpermeven}
  \pi=\biggl(\begin{smallmatrix}
    1 & 2 & 3 & \cdots & i    & \cdots & N/2 & N/2+1 & \cdots & j         & \cdots & N-2 & N-1 & N \\
    1 & 3 & 5 & \cdots & 2i-1 & \cdots & N-1 & N     & \cdots & 2N-2(j-1) & \cdots & 6   & 4   & 2
  \end{smallmatrix} \biggr),
\end{equation}
while for $N$ odd, the map $\pi$ is defined:
\begin{equation}\label{e:circpermodd}
  \pi=\biggl(\begin{smallmatrix}
    1 & 2 & 3 & \cdots & i  & \cdots & \left\lfloor N/2 \right \rfloor & \left \lceil N/2 \right \rceil & \cdots & j        & \cdots & N-2 & N-1 & N \\
    2 & 4 & 6 & \cdots & 2i & \cdots & N-1                             & N                              & \cdots & 2(N-j)+1 & \cdots & 5   & 3   & 1
  \end{smallmatrix}\biggr).
\end{equation}

It is clear that $\pi$ induces a digraph isomorphism $\varphi:\Lc\rightarrow\Dc$ which preserves edge colors. Moreover, $\varphi^{-1}:\Dc\rightarrow\Lc$ is induced by $\pi^{-1}\in S^{N}$.

%
%
The isomorphisms above enable us to state our main result, that $\L\cong\D$ as edge-colored digraphs.  
\begin{Thm}\label{t:main} $\Phi:\Lp\rightarrow\Dp$ defined by $\Phi=\gamtp\,\circ\,\gamct\,\circ\,\varphi\,\circ\,\gamtc\,\circ\,\gampt$ is an edge-color preserving isomorophism of directed graphs.  $\Phi^{-1}:\Dp\rightarrow\Lp$ is given by $\Phi^{-1}=\gamtp\,\circ\,\gamct\,\circ\,\varphi^{-1}\,\circ\,\gamtc\,\circ\,\gampt$.  
\end{Thm}
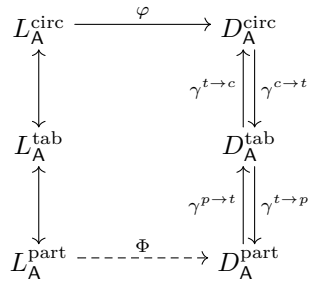
\begin{proof}\ 

\begin{center}
\begin{tikzcd}[row sep=3em,column sep=5em]
		\Lc
			\arrow[transform canvas={yshift=.5ex}]{r}{\varphi}
			\arrow[leftrightarrow]{d}
 & 	\Dc
			\arrow[transform canvas={xshift=.5ex}]{d}{\gamct}\\
		\Lt
			\arrow[leftrightarrow]{d}
 & 	\Dt
			\arrow[transform canvas={xshift=-.5ex}]{u}{\gamtc}
			\arrow[transform canvas={xshift=.5ex}]{d}{\gamtp}\\
		\Lp
			\arrow[dashrightarrow,transform canvas={yshift=.5ex}]{r}{\Phi}
 & \Dp
			\arrow[transform canvas={xshift=-.5ex}]{u}{\gampt}
\end{tikzcd}			

\end{center}

The isomorphisms on the left-hand side are clear from Definitions \ref{def:lpart}, \ref{def:ltab}, and \ref{def:lcirc}; as each of these are isomorphic to $L_{A}$.  
\end{proof}

The conclusion is that $\Da$, which is defined as an edge-colored directed game graph, actually has the structure of the diamond-colored distributive lattice $\La$.  This means that the four move-minimizing game graph questions set out in Section \ref{s:games} can all be answered in the affirmative.  Specifically, considering the non-colored version of the questions described in the paragraphs before Theorem \ref{t:MoveMinGameTheorem}, we have that $\D$ is connected; (1): There is an explicit formula for the distance between two domino partitions $\sigma$ and $\tau$ (questions (1) and (2) are combined in this setting); and $(3)$: An explicit path from $\sigma$ to $\tau$ can be prescribed.  Though we will focus on the non-colored situation below, the answer to the colored versions of questions (2) and (3) will be very clear with our method. See Examples \ref{ex:dominomoves} and \ref{ex:path}.

\begin{Ex}\label{ex:iso} Consider $(4,3)\in\D(2,4)$.  The diagram below shows the progression of this partition through the commutative diagram in Theorem \ref{t:main} to the associated partition $(1,1)\in\L(2,4)$.  

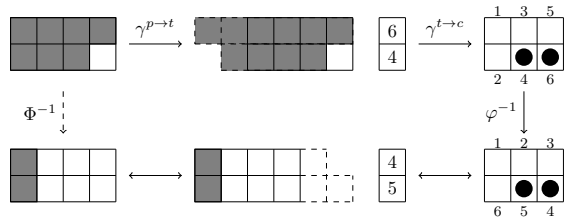
\begin{center}
\scalebox{.7}{%

\begin{tikzpicture}

	\filldraw[fill=gray] (-.5,0) rectangle (1,1);
	\filldraw[fill=gray] (1,.5) rectangle (1.5,1);

	\foreach \i in {0,...,2}
		\draw (-.5,\i/2) -- (1.5,\i/2);
	\foreach \j in {0,...,4}
		\draw (\j/2-.5,0) -- (\j/2-.5,1);
	
	\draw[->] (1.75,.5) -- node[above]{$\gampt$} (2.75,.5);
	
	\draw[dashed,->] (.5,-.4) -- node[left]{$\Phi^{-1}$} (.5,-1.2);
	
	\filldraw[fill=gray,dashed] (3,.5) rectangle (6,1);
	\filldraw[fill=gray,dashed] (3.5,0) rectangle (5.5,.5);

	\foreach \i in {0,...,2}
		\draw (4,\i/2) -- (6,\i/2);
	\foreach \j in {8,...,12}
		\draw (\j/2,0) -- (\j/2,1);r
	\draw[dashed] (3.5,.5) -- (3.5,1);

	\draw (6.5,0) -- (7,0);
	\draw (6.5,.5) -- (7,.5);
	\draw (6.5,1) -- (7,1);
	\draw (6.5,0) -- (6.5,1);
	\draw (7,0) -- (7,1);
		
	\node at (6.75,.25){4};
	\node at (6.75,.75){6};
	
	\draw[->] (7.25,.5) -- node[above]{$\gamtc$} (8.25,.5);
	
	\foreach \i in {0,...,2}
			\draw (8.5,\i/2) -- (10,\i/2);
		\foreach \j in {16,...,19}
			\draw (\j/2+.5,0) -- (\j/2+.5,1);

	\node[font=\footnotesize] at (8.75,1.12){1};
	\node[font=\footnotesize] at (9.25,1.12){3};
	\node[font=\footnotesize] at (9.75,1.12){5};

	\node[font=\footnotesize] at (8.75,-.17){2};
	\node[font=\footnotesize] at (9.25,-.17){4};
	\node[font=\footnotesize] at (9.75,-.17){6};
	
	\node[circle,fill=black] at (9.25,.25){};
	\node[circle,fill=black] at (9.75,.25){};

	\draw[->] (9.25,-.4) -- node[left]{$\varphi^{-1}$}(9.25,-1.2);
	
	\foreach \i in {0,...,2}
			\draw (8.5,\i/2-2.5) -- (10,\i/2-2.5);
		\foreach \j in {16,...,19}
			\draw (\j/2+.5,-2.5) -- (\j/2+.5,-1.5);

	\node[font=\footnotesize] at (8.75,-1.38){1};
	\node[font=\footnotesize] at (9.25,-1.38){2};
	\node[font=\footnotesize] at (9.75,-1.38){3};

	\node[font=\footnotesize] at (8.75,-2.67){6};
	\node[font=\footnotesize] at (9.25,-2.67){5};
	\node[font=\footnotesize] at (9.75,-2.67){4};

	\node[circle,fill=black] at (9.25,-2.25){};
	\node[circle,fill=black] at (9.75,-2.25){};
	
	\filldraw[fill=gray] (3,-2.5) rectangle (3.5,-1.5);

	\foreach \i in {0,...,2}
		\draw (3,\i/2-2.5) -- (5,\i/2-2.5);
	\foreach \j in {6,...,10}
		\draw (\j/2,-2.5) -- (\j/2,-1.5);
	\draw[dashed] (5,-2.5) -- (6,-2.5);
	\draw[dashed] (5,-2) -- (6,-2);
	\draw[dashed] (5,-1.5) -- (5.5,-1.5);
	\draw[dashed] (5.5,-2.5) -- (5.5,-1.5);
	\draw[dashed] (6,-2.5) -- (6,-2);

	\draw (6.5,-2.5) -- (7,-2.5);
	\draw (6.5,-2) -- (7,-2);
	\draw (6.5,-1.5) -- (7,-1.5);
	\draw (6.5,-2.5) -- (6.5,-1.5);
	\draw (7,-2.5) -- (7,-1.5);
		
	\node at (6.75,-1.75){4};
	\node at (6.75,-2.25){5};
	
	\draw[<->] (8.25,-2) -- (7.25,-2);

\filldraw[fill=gray] (-.5,-2.5) rectangle (0,-1.5);
	
	\foreach \i in {0,...,2}
		\draw (-.5,\i/2-2.5) -- (1.5,\i/2-2.5);
	\foreach \j in {0,...,4}
		\draw (\j/2-.5,-2.5) -- (\j/2-.5,-1.5);
	
	\draw[<->] (2.75,-2) -- (1.75,-2);
	
\end{tikzpicture}

}
\end{center}
\end{Ex}

\paragraph{Solution to the Domino Game.}  Suppose one is given two $k\times (N-k)$ partitions, and the idea is to obtain one from the other.  In $\L$, the moves needed to acquire one partition from another are obvious:  one simply adds or removes boxes from one until arriving at the other.  In $\D$, the domino moves may not be as clear.  So, if one is given two \emph{domino} partitions, by Theorem \ref{t:main}, a winning strategy in $\D$ is to actually play the game in $\L$ and perform the corresponding moves in $\D$. That is, given $\sigma$ and $\tau$ in $\D$, on simply adds or removes boxes to $\Phi^{-1}(\sigma)$ to obtain $\Phi^{-1}(\tau)$ in $\L$, and mirrors the moves in $\D$.  One may answer question (1) by determining the number of moves required in $\L$, and questions (2) and (3) by noting the color of each move insde $\L$. See Example \ref{ex:dominogame24}.

\begin{Ex}\label{ex:dominogame24} The diagram below shows the use of $\Phi$ in a successful version of the domino game given $\sigma=(4,3)$ and $\tau=(1,1)\in\D(2,4)$.  
\end{Ex}

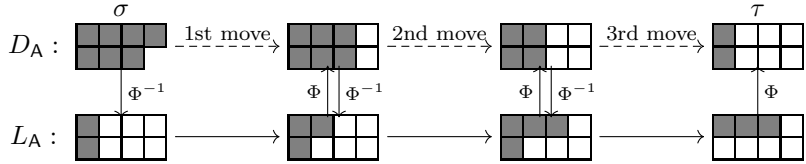
\begin{tikzcd}[row sep=2em,column sep=4em]
 &[-4em] \sigma & & & \tau \\[-2em] 
\D:&\ydiagram[*(gray)]{4,3}*[*(white)]{0,3+1}
		\arrow[dashed]{r}{\mbox{\footnotesize 1st move}}
		\arrow[]{d}{\Phi^{-1}}
& \ydiagram[*(gray)]{3,3}*[*(white)]{3+1,3+1}
		\arrow[dashed]{r}{\mbox{\footnotesize 2nd move}}
		\arrow[transform canvas={xshift=.5ex}]{d}{\Phi^{-1}}
& \ydiagram[*(gray)]{2,2}*[*(white)]{2+2,2+2}
		\arrow[dashed]{r}{\mbox{\footnotesize 3rd move}}
		\arrow[transform canvas={xshift=.5ex}]{d}{\Phi^{-1}}
& \ydiagram[*(gray)]{1,1}*[*(white)]{1+3,1+3}\\
\L:& \ydiagram[*(gray)]{1,1}*[*(white)]{1+3,1+3}
		\arrow{r}
&	\ydiagram[*(gray)]{2,1}*[*(white)]{2+2,1+3}
		\arrow{r}
		\arrow[transform canvas={xshift=-.5ex}]{u}{\Phi}
&	\ydiagram[*(gray)]{3,1}*[*(white)]{3+1,1+3}
		\arrow{r}
		\arrow[transform canvas={xshift=-.5ex}]{u}{\Phi}
&	\ydiagram[*(gray)]{3,0}*[*(white)]{3+1,4}
		\arrow{u}[swap]{\Phi}\\
\end{tikzcd}	

\paragraph{The diagonal coordinates and a linear algebraic approach.} We now turn our attention to the so-called ``diagonal coordinates'' of $\L$ and $\D$.  These seemingly awkward coordinates actually provide a very clean connection between the Domino Game and the classical distributive lattice $\L$, in that they help us prescribe answers to the non-colored versions of questions (1)--(3) above.  

As stated in Theorem \ref{t:ModularLatticeTheorem}, if one knows the rank of an element in a topographically balanced distributive lattice, then calculating the length of the shortest path between two elements is straightforward.  In $\L$, the rank of a given partition is simple:  It is the sum of the rows of the partition.  That is, if $\sigma=\sigma_1,\sigma_2,\ldots,\sigma_k$, $\rho(\sigma)=\sum \sigma_i$.  

The calculation of the corresponding rank of a partition $\tau$ in $\D$ is not so obvious.  Now, it is clear from Theorem \ref{t:main} that one can figure the rank of $\tau\in\D$ by computing the rank of $\Phi^{-1}(\tau)\in\L$, but it's not clear how to calculate the rank of an element without ``leaving'' $\D$ and relying on the classical object $\L$.  In terms of the domino game, this means that we cannot directly find the answer to (1)--(3) without utilizing $\Phi$ in the manner of Example \ref{ex:dominogame24}.  The diagonal coordinates, however, give the key to answering (1)--(3) without moving back and forth between $\L$ and $\D$.  

Recall again that the diagonal coordinates of both $\L$ and $\D$ represent each partition as vector in $\BZ^{N-1}$.  It is clear that the set $\varepsilon_{i}\in\BZ^{N-1}$, for $i \in[N-1]$ forms a $\BZ$-basis for the elements of $\Ld$.  That is, each partition $\sigma\in\L$, understood via it's diagonal coordinates, $\ds$, can be viewed uniquely as a linear combination of the basis described above, i.e. $\ds=\sum c_i\varepsilon_i$.  Moreover, the basis vectors $\left\{\varepsilon_i\right\}$ correspond to edges in $\L$, or, in our game graph terminology, `moves' in $\L$.  This means that $\sigma$, when viewed through its diagonal coordinates $\ds$, is a linear combination of the basis vectors that encode the colored `moves' used to reach $\sigma$ from the empty, minimum partition.  So, since the rank of a partition is simply the minimum number of moves required to go from the minimum partition to the given partition, the rank of $\sigma$ is given by $\rho(\sigma)=\sum c_i$.  

\begin{Ex}\label{ex:1}  Let $\sigma=(3,3)\in \L(2,4)$.  As in Example \ref{ex:big24}, we see that the diagonalization of $\sigma$ is $\ds =(0,1,2,2,1)$, and, according to the same picture, we see that 
\begin{align}\label{e:moves}
(0,1,2,2,1)&=\varepsilon_2+2\varepsilon_3+2\varepsilon_4+\varepsilon_5,
\end{align}
and has rank 6, since there are six moves from the minimum up to $\sigma$, one each of colors 2 and 5, two each of colors 3 and 4, regardless of the path traveled from the empty partition.
\end{Ex}

\paragraph{Diagonal moves.}  Under $\Phi:\L\rightarrow\D$, $\varepsilon_i$ corresponds to $\beta_i$, for each $i$, where the $\beta_i$'s, defined in Defintion \ref{def:ddiag} correspond to specific-colored the domino moves.  This means that for $\tau\in\D$ with diagonal coordinates $\dt$, we can encode the domino moves needed to go from the minimum partition in $\D$ to $\tau$ as a linear combination of $\beta_i$'s.  In particular, $\dt=\sum d_i\beta_i + \m$ where $\m$ denotes the diagonal coordinates of the minimum element of $\D$.  Moreover, we can calculate the rank of $\tau$ by summing the coefficients of the linear description of $\dt-\m$ while the $d_i$'s determine the minimum number of each colored required in the path from the minimum to $\tau$.  See Example \ref{ex:2}.

\begin{Ex}\label{ex:2}  As in Example \ref{ex:1}, but now considering domino coordinates, we have that $(0,1,2,2,1)\xrightarrow{\Phi}(0,1,1,2,1)$.  And, equivalently to Equation \ref{e:moves}, one can confirm with $D(2,4)$ in Example \ref{ex:big24} that 
\begin{align}\label{e:dommoves}
\begin{split}
(0,1,1,2,1)=&(0,-1,-1,0,0)+2(-1,0,0,0,0)+2(1,1,0,0,0)+(0,0,1,1,0)\\
 &+\left[(0,0,1,1,1)\right]\\
           =&\beta_2+2\beta_3+2\beta_4+\beta_5+\m
\end{split}
\end{align}
which again yields 6 moves form the minimum to $\tau$, with one each of colors 2 and 5, two each of colors 3 and 4, independent of the path selected starting at $\m$ in $\D(2,4)$.
\end{Ex}

\paragraph{The matrix of $\Phi$.}  The above means that the map $\Phi:\L\rightarrow\D$ can be described with an $(N-1)\times(N-1)$ matrix $P$ whose columns are the vectors $\beta_i$, $i\in[N-1]$, which vary according to $N$ even or odd.  But, because the 0-vector (empty partition) is not the minimum element in $\D$, $P$ must be adjusted.  That is, for $\ds\in\Ld$ corresponding to the partition $\sigma$, $P(\ds)+\m\in\Dd$ is the image of $\ds$ corresponding to the partition $\Phi(\sigma)$, where $\m$ is the minimum element of $\D$ in diagonal coordinates.  Let $\left[P\right]$ be the map induced by the matrix $P$, including the appropriate shift.  Then the commutative diagram in Theorem \ref{t:main} is extended to include the following: 
\begin{center}
\begin{tikzcd}[row sep=3em,column sep=5em]
	\Lp
			\arrow[transform canvas={yshift=.5ex}]{r}{\Phi}
			\arrow[leftrightarrow]{d}
 & \Dp
			\arrow[transform canvas={xshift=.5ex}]{d}{\gampd}\\
		\Ld
			\arrow[dashrightarrow,transform canvas={yshift=.5ex}]{r}{\left[P\right]}
 &	\Dd
			\arrow[transform canvas={xshift=-.5ex}]{u}{\gamdp}
\end{tikzcd}			

\end{center}
The matrix $P$ is instrumental to playing the domino game without `leaving' the domino environment.  Indeed, given $\sigma\in\D$ with diagonalization $\ds$, it may not be clear how to see $\ds$ as a sum of the $\beta_i$'s.  But, since the minimum $\m$ of $\D$ is understood (see Definition \ref{def:ddiag}), we can find $(\ds-\m)$.  Then $P^{-1}(\ds-\m)=\dt\in\BZ^{N-1}$ gives the $\BZ$-linear combination of $\beta_i$'s for $(\ds-\m)$; for then $(\ds-\m)=\sum_i c_i\beta_i$.  It is this linear combination that can be used to find answers to (1)-(2) and (3).  See Example \ref{ex:matrixmult}.

\begin{Ex}\label{ex:matrixmult} We now connect examples \ref{ex:1} and \ref{ex:2}.  Considering $D(2,4)$, let $\ds=(0,1,1,2,1)\in\Dd$ as in Equation \ref{e:dommoves}, with $\m=(0,0,1,1,1)$.  Then
\[\begin{array}{cc}
	P^{-1}(\ds-\m) = & \left(\begin{array}{ccccc}
												0 & 0 & -1 & 1 & 0 \\
												0 & -1 & 0 & 1 & 0 \\
												0 & -1 & 0 & 0 & 1 \\
												-1 & 0 & 0 & 0 & 1 \\
												-1 & 0 & 0 & 0 & 0 \end{array}\right)^{-1}\left(\begin{array}{c} 0 \\ 1 \\ 0 \\ 1 \\ 0\end{array}\right)=\left(\begin{array}{c} 0 \\ 1 \\ 2 \\ 2 \\ 1\end{array}\right)\end{array}.\]
So, as in Equation \ref{e:dommoves}, we see that $\ds=\beta_2+2\beta_3+2\beta_4+\beta_5$.
\end{Ex}

Note that, because $[P]:\Ld\rightarrow\Dd$, we do technically `leave' the domino environment for the simpler $\L$ environment, but this matrix encodes all of the information to answer the questions we have posed without explicitly calculating the corresponding partitions in $\L$.  

\paragraph{A Domino solution to (1)--(3).}  Given partitions $\sigma$ and $\tau$ in $\D$ with corresponding diagonal coordinates $\ds$ and $\dt$, and with $\m$ be the minimum element in $\D$, consider the expansions $(\ds-\m)=\sum c_i\beta_i$ and $(\dt-\m)=\sum d_i\beta_i$ which can be found by $P^{-1}(\ds-\m)$ and $P^{-1}(\dt-\m)$.  We use the mountain path defined in Section \ref{s:lattices} to determine the length of the shortest path from $\sigma$ to $\tau$ and we prescribe that path.  Let 
\begin{align*}
\mathsf{S}&=\{i^{c_i}\mid 1\leq i\leq N-1\}, \text{ and }\\
\mathsf{T}&=\{i^{d_i}\mid 1\leq i\leq N-1\}
\end{align*}
be the multisets of the vector expansions of $\ds$ and $\dt$, where $i$ has multiplicity $c_i$ and $d_i$ in each set, respectively.  These multisets encapsulate the moves required to get from the minimum to $\sigma$ and $\tau$.  Let $\mathsf{u}$ be the diagonal coordinates of $\sigma\vee\tau$.  It is clear that $(\mathsf{u}-\m)$ has multiset $\mathsf{S}\cup\mathsf{T}$.  Then since one can get from $\sigma$ to $\tau$ in a minimal number of moves by traveling through the join (the mountain path), we have that the number of moves required to get from $\sigma$ to $\tau$ is 
\begin{equation}\label{e:countdominomoves}
|(\mathsf{S}\cup\mathsf{T})-\mathsf{S}|+|(\mathsf{S}\cup\mathsf{T})-\mathsf{T}|
\end{equation}
while the required number of $i$-colored moves is $c_i+d_i$, for $i \in [N-1]$.

\begin{Ex}\label{ex:dominomoves} Let $\sigma=(4,4)$ and $\tau=(1,1)$ in $\mathcal{D}_A(2,4)$.  Then $\ds=(1,2,2,2,1)$, $\dt=(0,0,0,1,1)$, $\mathsf{S}=\{3,4,4,5\}$ and $\mathsf{T}=\{2,3,4\}$. (Refer to Figure \ref{fig:D24Full}.) $\sigma\vee\tau$ has associated multiset $\{2,3,4,4,5\}$ and we see  that one can get from $\sigma$ to $\tau$ in three moves.  And in Equation \ref{e:countdominomoves}, we have 
\[|(\mathsf{S}\cup\mathsf{T})-\mathsf{S}|+|(\mathsf{S}\cup\mathsf{T})-\mathsf{T}|=1+2=3\]
\end{Ex}

We can also prescribe a shortest path from $\sigma$ to $\tau$ that uses a minimum number of $i$-colored edges, for each $i$.  Indeed, starting from $\sigma$, we describe the procedure in terms of steps.
\begin{enumerate}
\setcounter{enumi}{-1}
	\item Calculate $P^{-1}(\ds-\m)$ and $P^{-1}(\dt-\m)$ to find the multisets $\mathsf{S}$ and $\mathsf{T}$ as above.
	\item Put $\ds^1=\ds+\beta_{i_1}$ where $i_1$ is the minimal element of $((\mathsf{S}\cup\mathsf{T})-\mathsf{S})$ for which $\ds+\beta_{i_1}$ is the diagonalization of a partition.  (That is, $\beta_{i_1}$ corresponds to a legal domino move.) 
	\item Put $\ds^2=\ds^1+\beta_{i_2}$ where $i_2$ is the minimal element of $((\mathsf{S}\cup\mathsf{T})-\mathsf{S})-\{i_1\}$ for which $\ds^1+\beta_{i_2}$ is the diagonalization of a partition.  
	\item Let $\ds^3=\ds^2+\beta_{i_3} \dots$

\hspace*{.5cm}\vdots
\end{enumerate}
This process must terminate at $\sigma\vee\tau$, since its corresponding multiset is $\mathsf{S}\cup\mathsf{T}$.  The chosen sequence $\left\{\beta_{i_l}\right\}_{l=1}^n$ then defines an edge path from $\sigma$ to $\sigma\vee\tau$.  Repeat the process using $(\mathsf{S}\cup\mathsf{T})-\mathsf{T}$ to obtain a sequence $\left\{\beta_{j_k}\right\}_{k=1}^m$ and a path from $\tau$ to $\sigma\vee\tau$.  Then the first path followed by the second in reverse defines a path from $\sigma$ to $\tau$:
	\begin{enumerate}
	\setcounter{enumi}{-1}
	\item $\sigma\longleftrightarrow \ds$
	\item $\ds^1=\ds+\beta_{i_1}$
	
	\hspace*{.5cm}\vdots
	
	\item[$n$.] $\ds^n=\ds^{n-1}+\beta_{i_n}=\dt^m\longleftrightarrow\sigma\vee\tau$
	\item[$n+1$] $\dt^{m-1}=\dt^m-\beta_{j_m}$
	
	\hspace*{.5cm}\vdots
	
	\item[$n+m.$] $\dt=\dt^1-\beta_{j_1}\longleftrightarrow\tau$
	\end{enumerate}
Note that this is the mountain path described in Section \ref{s:lattices}; and this type of strategy was executed in Example \ref{ex:dominogame24}.

\begin{Ex}\label{ex:path} Let $\sigma$ and $\tau$ be as in Example \ref{ex:dominomoves}.  Then the path from $\sigma$ to $\sigma\vee\tau$ consists of the edge $2$.  The path from $\tau$ to $\sigma\vee\tau$ consists of traversing edge $5$ then edge $4$.  Therefore the path beginning at $\sigma$ consisting of edges $2$, $4$, and then $5$ in gives a path from $\sigma$ to $\tau$.  That is, 
\begin{enumerate}
\setcounter{enumi}{-1}
\item $\sigma=(4,4)\longleftrightarrow \ds=(1,2,2,2,1)$
\item $\ds^1=\ds+\beta_2=(1,1,1,2,1)=\dt^2\longleftrightarrow\sigma\vee\tau$
\item $\dt^1=\dt^2-\beta_4=(0,0,1,2,1)$
\item $\dt=\dt^1-\beta_5=(0,0,0,1,1)\longleftrightarrow\tau=(1,1)$.
\end{enumerate}
\end{Ex}

\bibliography{latticebib}

\end{document}